\documentclass[preprint,11pt]{elsarticle}
\usepackage{mathrsfs}
\usepackage{amsfonts}
\usepackage{amssymb}
\usepackage{amsmath}
\usepackage{color}
\usepackage{graphicx}
\usepackage{galois}

\numberwithin{equation}{section}

\newtheorem{theorem}{Theorem}[section]
\newtheorem{lemma}{Lemma}[section]
\newtheorem{proposition}{Proposition}[section]

\newtheorem{definition}{Definition}[section]
\newproof{pf}{Proof}
\newtheorem{remark}{Remark}[section]

\begin{document}

\begin{frontmatter}

\title{Existence of Entropy Solutions to Two-Dimensional Steady Exothermically Reacting Euler Equations}
\author{Gui-Qiang Chen, Changguo Xiao \& Yongqian Zhang}
\address{Mathematical Institute, University of Oxford,
Oxford, OX2 6GG, UK;\\
 School of Mathematical Sciences, Fudan University,
 Shanghai 200433, China.\\
 E-mail: chengq@maths.ox.ac.uk}
\address{
         School of Mathematical Sciences, Fudan University, Shanghai 200433, China. \\
         E-mail: 09110180021@fudan.edu.cn; yongqianz@fudan.edu.cn}

\begin{abstract}
We are concerned with the global existence of entropy solutions of the two-dimensional steady Euler
equations for an ideal gas, which undergoes a one-step exothermic chemical reaction
under the Arrhenius-type kinetics.
The reaction rate function $\phi(T)$ is assumed to have a positive lower bound.
We first consider the Cauchy problem (the initial value problem), that is,
seek a supersonic downstream reacting flow when the incoming flow is supersonic,
and establish the global existence of entropy solutions
when the total variation of the initial data is sufficiently small.
Then we analyze the problem of steady supersonic, exothermically reacting Euler flow past a Lipschitz wedge,
generating an additional detonation wave attached to the wedge vertex, which can be then formulated
as an initial-boundary value problem.
We establish the globally existence of entropy solutions containing the additional
detonation wave (weak or strong, determined by the wedge angle at the wedge vertex)
when the total variation of both the slope of the wedge boundary and the incoming
flow is suitably small.
The downstream asymptotic behavior of the global solutions is also obtained.


\end{abstract}

\end{frontmatter}

\bigskip
\section{Introduction}
We are concerned with
the two-dimensional steady supersonic Euler flow of an exothermically
reacting ideal gas, which is governed by
\begin{eqnarray}
(\rho u)_x+(\rho v)_y = 0,\label{d1}\\
(\rho u^2+p)_x+(\rho uv)_y = 0,\\
(\rho uv)_x+(\rho v^2+p)_y = 0,\\
\big((\rho E+p)u\big)_x+\big((\rho E+p)v\big)_y = 0,\\
(\rho uZ)_x+(\rho vZ)_y = -\rho Z\phi(T). \label{d2}
\end{eqnarray}
Here $(u,v)$ is the velocity, $p$ the scalar pressure, $\rho$ the density, $Z$ the fraction of unburned
gas in the mixture, $\phi(T)$ the reaction rate, $q$ the specific binding energy of unburned gas,
and $E=\frac{1}{2}(u^2+v^2)+e(\rho,p)+qZ$ the specific total energy with the specific internal energy $e$
that is a given function of $(\rho,p)$
defined through thermodynamical relations.

For an ideal gas,
\begin{equation}
p=R\rho T,\quad e=c_vT,\quad \gamma=1+\frac {R}{c_v}>1,
\end{equation}
where $R$ and $c_v$ are positive constants, and $\gamma$ is the adiabatic exponent.
We identify $c_v+R=c_p$ as the specific heat at constant pressure.

\par We assume for simplicity that the specific heats and molecular weights of the reactant
and product gases are the same and that the reaction rate function $\phi$ is monotonically increasing
and Lipschitz continuous. In addition, inadmissible discontinuous solutions are eliminated
by requiring the following entropy condition:
\begin{equation}
(\rho uS)_x+(\rho vS)_y \ge \frac{q\rho Z \phi(T)}{T}.
\end{equation}

\par We first consider the Cauchy problem (the initial value problem) for \eqref{d1}--\eqref{d2}
in the region $\{x\ge 0, y\in \mathbb{R}\}$, with initial incoming flow (initial data):
\begin{equation}
(u,v,p,\rho,Z)(0,y)=(u_0,v_0,p_0,\rho_0,Z_0)(y),\qquad \text{$y \in \mathbb{R}$}.\label{d3}
\end{equation}
We assume that $u_0(y),v_0(y),p_0(y),\rho_0(y)$, and $Z_0(y)$ are bounded and have bounded total
variation with $Z_0(-\infty)=\lim_{y \to -\infty}Z_0(y)=0$.
We further assume that there are positive constants $u', \rho '$, and $T'$ such that
\begin{equation}\label{d3-a}
u_0 > c_0\ge u'>0,\qquad \rho_0\ge \rho '>0,\qquad T_0\ge T'>0,
\end{equation}
where  $c=\sqrt{\frac{\gamma p}{\rho}}$ is the local sonic speed.
We make this assumption on the initial data to ensure that the flow
is supersonic ({\it i.e.} $u^2+v^2>c^2$).
\bigskip

\begin{center}
\setlength{\unitlength}{1mm}
\begin{picture}(60,60)(0,-30)
\linethickness{1pt}
\put(0,0){\vector(1,0){60}}
\put(0,-30){\vector(0,1){60}}
\thinlines
\put(0,8){\line(1,1){20}}
\put(0,-4){\line(1,1){20}}
\put(0,-16){\line(1,1){20}}
\put(0,-28){\line(1,1){20}}
\thinlines
\put(-4,-1){$O$}
\put(61,-2){$x$}
\put(-3,28){$y$}
\put(-15,14){$u_0>c_0$}
\put(-15,8){$v_0$}
\put(-15,2){$p_0$}
\put(-15,-4){$\rho_0$}
\put(-15,-10){$Z_0$}
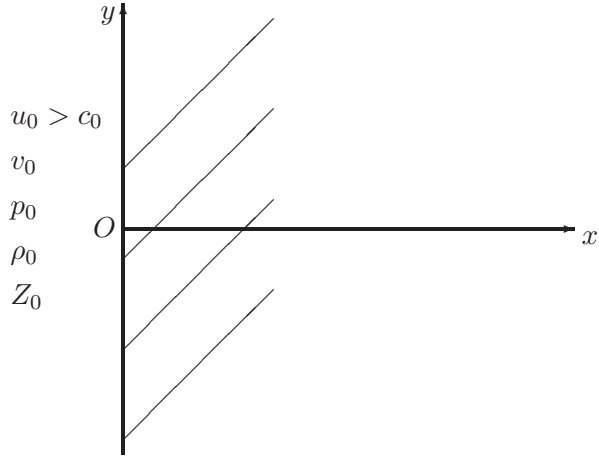
\put(-20,-34){Fig.1.  Supersonic Euler flow through the left boundary $x=0$ }
\end{picture}
\end{center}

\bigskip
\medskip

We assume the initial data to be such that the reaction rate function $\phi(T)$ never vanishes,
so that there is a positive minimum value $\Phi:=\phi(T^{\prime})>0$.
In a sense, this is a very realistic condition. Typically, $\phi (T)$ has the Arrhenius form:
\begin{equation}
\phi(T)=T^{\alpha}e^{-\frac{E}{RT}},
\label{arr}
\end{equation}
which vanishes only at absolutely zero temperature, where $\alpha$ is a positive constant.
We make this assumption in order to obtain the uniform decay of the reactant to zero.
Although the total variation of the solution may very well increase while the reaction is active,
the reaction must eventually die out along the flow trajectories.
Consequently, the increase in total variation can be estimated rigorously.

\par In Chen-Wagner \cite{chen1},
the large-time existence of entropy solutions to the Cauchy problem has been established for
the time-dependent equations of planar flow of an exothermically reacting ideal gas.
%
%
The total variation of the initial data is bounded by a parameter
$\epsilon=\gamma-1$, which grows arbitrarily large as $\epsilon\to 0$ whose limiting case
is the isothermal gas.
Global entropy solutions are obtained
by using the Glimm fractional-step scheme based on the Glimm scheme.

In this paper, we first
establish a global existence theory for entropy solutions of the Cauchy problem
for two-dimensional, exothermically
reacting steady Euler equations by further developing the Glimm scheme, under the condition
that the total variation of the initial data in \eqref{d3} is small.
Then this approach is further developed for solving
the supersonic reacting Euler flow past Lipschitz wedges.
For a non-reacting supersonic flow past a straight wedge, an attached plane shock is generated at
the wedge vertex. When the supersonic flow is governed by the exothermically reacting
steady Euler equations,
the attached detonation wave is no longer a plane wave even for the straight wedge,
whose strength (weak or strong) is determined by the wedge angle and the incoming flow.
Nevertheless, we establish that, when the total variation of both the incoming supersonic flow
and the slope of the wedge boundary is suitably small, there exists a global entropy solution
containing the (weak or strong) detonation wave.
The downstream asymptotic behavior of entropy solutions is also obtained.

The organization of this paper is as follows.
In Section 2, we discuss some basic features of the exothermically reacting Euler equations
\eqref{d1}--\eqref{d2}.
The Glimm fractional-step scheme is described for the Cauchy problem \eqref{d3}
for system \eqref{d1}--\eqref{d2} in Section 3.
In Section 4, we establish uniform bounds on the total variation in the $y$-direction
of the Glimm fractional-step approximate solutions for the Cauchy problem \eqref{d3}.
In Section 5, we establish uniform bounds on the total variation
of the Glimm fractional-step approximate solutions in the $y$--variable for the
initial-boundary value problem  \eqref{dd3}--\eqref{dd4} for \eqref{d1}--\eqref{d2}
concerning the supersonic reacting Euler flow past Lipschitz wedges,
when the wedge angle at the wedge vertex
is small.
In Section 6, the convergence of approximate solutions to an entropy solution is established
for both the Cauchy problem
\eqref{d3} and the initial-boundary value problem  \eqref{dd3}--\eqref{dd4}
for \eqref{d1}--\eqref{d2}.
The downstream asymptotic behavior of entropy solutions is also clarified in Section 7.
In Section 8, we extend the results in Sections 5--6 for the case of small wedge angle to
the case of large wedge angle, for which the entropy solution contains a strong detonation wave
generated between the incoming fluid and the wedge boundary at the wedge vertex.

\section{Basic Features of the Exothermically Reacting Euler Equations}

In this section, we discuss some basic features of system  \eqref{d1}--\eqref{d2}.

\subsection{\textbf{Euler equations}}

\par System \eqref{d1}--\eqref{d2} can be rewritten in the following form:
\begin{equation}
W(U)_x+H(U)_y=G(U),
\label{2}
\end{equation}
with $U=(u,v,p,\rho,Z)$, where
\begin{eqnarray*}
&& W(U)=(\rho u,\rho u^2+p,\rho uv,\rho u(\bar{h}+\frac{u^2+v^2}{2}),\rho uZ),\\[1.5mm]
&& H(U)=(\rho v,\rho uv,\rho v^2+p,\rho v(\bar{h}+\frac{u^2+v^2}{2}),\rho vZ),\\[1.5mm]
&& G(U)=(0,0,0, q\rho \phi(T)Z,-\rho \phi(T)Z)
\end{eqnarray*}
with $\bar{h}=\frac {\gamma p}{(\gamma -1)\rho}$.

\par In the case when $G(U)$ is identically zero, system \eqref{2} becomes
a system of conservation laws:
\begin{equation}
W(U)_x+H(U)_y=0.
\label{3}
\end{equation}
For a smooth solution $U(x,y)$, system \eqref{3} is equivalent to
\begin{equation}
\nabla_UW(U)U_x+\nabla_UH(U)U_y=0.
\end{equation}
Then the eigenvalues of \eqref{3} are the roots of the 5th order polynomial
\begin{equation}
\mathrm{det}(\lambda\nabla_UW(U)-\nabla_UH(U)),
\end{equation}
that is, the solutions of the equation:
\begin{equation}
(v-\lambda u)^3\big((v-\lambda u)^2-c^2(1+\lambda ^2)\big)=0,
\end{equation}
where $c=\sqrt{\frac{\gamma p}{\rho}}$ is the sonic speed.

If the flow is supersonic ({\it i.e.} $u^2+v^2>c^2$),
system \eqref{2} is hyperbolic.
In particular, when $u>c$, the system has five eigenvalues
in the $x$-direction:
\begin{equation}
\lambda _i=\frac{v}{u}, \quad i=2,3,4;
\qquad\quad\, \lambda_j=\frac{uv+(-1)^{\frac{j+3}{4}} c\sqrt{u^2+v^2-c^2}}{u^2-c^2},
\quad j=1,5,
\end{equation}
and the corresponding linearly independent eigenvectors:
\begin{equation}
r_j=\kappa_j(-\lambda_j,1,\rho(\lambda_j u-v), \frac{\rho (\lambda_ju-v)}{c^2},0)^\top, \quad j=1,5;
\end{equation}
\begin{equation}
r_2=(u,v,0,0,0)^\top, \quad r_3=(0,0,0,\rho,0)^\top, \quad r_4=(0,0,0,0,\frac{1}{\rho u})^\top,
\end{equation}
where $\kappa_j$ are chosen so that $r_j\cdot\nabla \lambda_j=1$ since the $j$th-characteristic
fields are genuinely nonlinear, $j=1,5$.
Note that $r_j\cdot\nabla \lambda _j=0,j=2,3,4$, that is, these characteristic fields are always
linearly degenerate.

In particular,
at a constant state $\tilde{U}=(\tilde{u},0,\tilde{p},\tilde{\rho}, \tilde{Z})$,
$$
\lambda_2(\tilde{U})=\lambda_3(\tilde{U})=\lambda_4(\tilde{U})=0,
\qquad \lambda_1(\tilde{U})=-\frac{\tilde{c}}{\sqrt{\tilde{u}^2-\tilde{c}^2}}=-\lambda_5(\tilde{U})<0.
$$

\begin{definition}[Entropy Solutions]
A function $U=U(x,y)\in BV(\mathbb{R}^+\times \mathbb{R})$ is called an entropy solution of
problem \eqref{d3} for system \eqref{d1}--\eqref{d2} provided that
\begin{enumerate}
\item[\rm (i)] $U$ is a weak solution of problem \eqref{d3} for system \eqref{d1}--\eqref{d2}, that is,
\begin{equation}
\int_{-\infty}^{\infty}\int_{0}^{\infty}
\big(W(U)\phi_x+H(U)\phi_y+G(U)\phi\big)\, dxdy+\int_{-\infty}^{\infty}W(U_0(y))\phi(0,y)\, dy= 0
\end{equation}
for any $\phi \in C_0^{\infty}([0,\infty)\times(-\infty,\infty))$;

\medskip
\item[\rm (ii)] {For any convex entropy pair $(\eta,q)$ with respect to $W(U)$, the following inequality}
\begin{equation}
\eta(W(U))_x+q(W(U))_y\le \nabla_W\eta(W(U))G(U)
\end{equation}
holds in the sense of distributions, that is,
\begin{eqnarray}
\int_{-\infty}^{\infty}\int_{0}^{\infty}\big(\eta(W(U))\phi_x+q(W(U))\phi_y+\nabla_W \eta(W(U))G(U)\phi\big)\,dxdy\\
+\int_{-\infty}^{\infty}\eta(W(U_0(y)))\phi(0,y)\, dy\ge 0
\end{eqnarray}
for any $\phi \in C_0^{\infty}([0,\infty)\times(-\infty,\infty))$ and $\phi(x,y)\ge 0$.
\end{enumerate}
\end{definition}

\begin{remark}
In particular, $\eta(W)=-\rho uS$ is an entropy which is convex with respect to $W$,
while $q(W)=-\rho vS$ is the corresponding entropy flux, when $u>c>0$.
\end{remark}

\par As in \cite{chen1}, if we rewrite system \eqref{3} in Lagrangian coordinates:
\begin{eqnarray}
(x', m)=(x, m(x,y))
\end{eqnarray}
with
$\mathrm{d}m=\rho u \mathrm{d}y-\rho v \mathrm{d}x$,
then the fifth equation in \eqref{3} becomes
\begin{equation}
Z_{x'}=0.
\end{equation}
It states that the $Z$-component is decoupled from $(u,v,p,\rho)^\top$
in the solution of the non-reacting Riemann problem.

\subsection{\textbf{Wave curves in the phase space}}
We now analyze some basic properties of nonlinear waves.
We focus on the case when $u>c>0$ in the state space.
Seek the self-similar solutions to system \eqref{3}:
\begin{equation}
(u,v,p,\rho,Z)(x,y)=(u,v,p, \rho,Z)(\xi),\quad \xi =\frac {y}{x},
\end{equation}
which connect to a fixed constant state $U_0=(u_0,v_0,p_0,\rho_0,z_0)$.
Then we have
\begin{equation*}
\mathrm{det}\big(\xi \nabla_UW(U)-\nabla_UH(U)\big)=0,
\label{4}
\end{equation*}
which implies
\begin{equation}
\xi =\lambda_i(U)=\frac{v}{u}, \,\,\, i=2,3,4;\quad \text{or}\quad
\xi =\lambda_j(U),
\,\,\, j=1,5.
\end{equation}

\par Plugging $\xi=\lambda_i(U), i=2,3,4$, into \eqref{4}, we obtain
\begin{equation*}
dp=0,\qquad  vdu-udv=0,
\end{equation*}
which yields the contact discontinuity curves $C_i(U_0)$ in the phase space:
\begin{equation*}
C_i(U_0):\,\, p=p_0,\, w=\frac{v}{u}=\frac{v_0}{u_0}, \qquad i=2,3,4.
\end{equation*}
More precisely, we have
\begin{equation}
C_2(U_0):\,\, U=(u_0e^{\sigma_2},v_0e^{\sigma_2},p_0,\rho_0,Z_0)^\top \label{a},
\end{equation}
with strength $\sigma_2$ and slope $\frac{v_0}{u_0}$, which is determined by
 \begin{equation}
 \left\{
 \begin{array}{ll}
  \frac{dU}{d\sigma_2}=r_2(U),\\ [2mm]
   U|_{\sigma_2=0}=U_0;
 \end{array} \right.
 \end{equation}
 and
\begin{equation}
C_3(U_0):\,\, U=(u_0,v_0,p_0,\rho_0e^{\sigma_3},Z_0)^\top,\label{b}
\end{equation}
with strength $\sigma_3$ and slope $\frac{v_0}{u_0}$, which is determined by
\begin{equation}
\left\{
 \begin{array}{ll}
  \frac{dU}{d\sigma_3}=r_3(U),\\[2mm]
   U|_{\sigma_3=0}=U_0;
 \end{array} \right.
 \end{equation}
 and
\begin{equation}
C_4(U_0):\,\, U=(u_0,v_0,p_0,\rho_0,Z_0+\frac{\sigma_4}{\rho_0u_0})^\top,
\label{c}
\end{equation}
 with strength $\sigma_4$ and slope $\frac{v_0}{u_0}$, which is determined by
\begin{equation}
\left\{
 \begin{array}{ll}
  \frac{dU}{d\sigma_4}=r_4(U),\\[2mm]
   U|_{\sigma_4=0}=U_0.
 \end{array} \right.
 \end{equation}
We can see that $\sigma_4$ is the difference between
$w_5=\rho uZ$ in the Riemann problem.

\par Plugging $\xi=\lambda_j(U), j=1,5$, into \eqref{4}, we obtain the $j$-th rarefaction
wave curve $R_j(U_0)$, $j=1,5$,  in the phase space through $U_0$:
 \begin{equation}
 R_j(U_0):\,\, dp=c^2d\rho,\,\, du=-\lambda_jdv,\,\, \rho(\lambda_ju-v)dv=dp,\,\, dZ=0,\qquad  j=1,5.
 \end{equation}

\par Now we consider discontinuous solutions so that the equations in \eqref{3}
are satisfied in the distributional sense.
This implies that the following Rankine-Hugoniot conditions hold along the discontinuity
with speed $s$, which connects
to a state $U_0=(u_0,v_0,p_0,\rho_0,Z_0)$:
 \begin{alignat}{2}
 s[\rho u]&=[\rho v], \label{5}\\[1.5mm]
 s[\rho u^2+p]&=[\rho uv], \\[1.5mm]
 s[\rho uv]&=[\rho v^2+p], \\[1.5mm]
 s[\rho u(\bar h+\frac{u^2+v^2}{2})]&=[\rho v(\bar h+\frac{u^2+v^2}{2})],\\[1.5mm]
 s[\rho uZ]&=[\rho vZ], \label{6}
 \end{alignat}
where the jump symbol $[\cdot]$ stands for the value of the quantity of the front-state minus that of the back-state.
Then we have
\begin{equation*}
 (v_0-su_0)^3\big((v_0-su_0)^2-\bar{c}^2(1+s^2)\big)=0,
\end{equation*}
where $\bar{c}^2=\frac{c_0^2}{b}\frac{\rho}{\rho_0}$ and $b=\frac{\gamma+1}{2}-\frac{\gamma-1}{2}\frac{\rho}{\rho_0}$.
This implies
\begin{equation}
s=s_i=\frac{v_0}{u_0}, \qquad i=2,3,4,
 \end{equation}
 or
\begin{equation}
s=s_j=\frac{u_0v_0+(-1)^{\frac{j+3}{4}}\bar{c}\sqrt{u_0^2+v_0^2-\bar{c}^2}}{u_0^2-\bar{c}^2},\qquad j=1,5,
\end{equation}
where $u_0>\bar{c}$ for small shocks.

\par Plugging $s_i$, $i=2,3,4$, into \eqref{5}--\eqref{6},
we obtain the same $C_i(U_0)$, $i=2,3,4$, as defined in \eqref{a}, \eqref{b}, and \eqref{c};
while plugging $s_j$, $j=1,5$, into \eqref{5}--\eqref{6},
we obtain the $j$th shock wave curve $S_j(U_0)$, $j=1,5$, through $U_0$:
\begin{equation}
S_j(U_0):\, \, [p]=\frac{c_0^2}{b}[\rho],\,\, [u]=-s_j[v], \,\, \rho_0(s_ju_0-v_0)[v]=[p],\,\, [Z]=0.
\end{equation}
Note that the shock wave curve $S_j(U_0)$ contacts with $R_j(U_0)$ at $U_0$ up to second order.

\par Following Lax \cite{lax}, we can parameterize any physically admissible wave curve
in a neighborhood of a constant $\tilde{U}$, $O_{\epsilon}(\tilde{U})$,
by $\sigma_j \mapsto \Phi_j(\sigma_j,U_b)$, with $\Phi_j \in C^2, \Phi_j|_{\sigma_j=0}=U_b$,
and $\frac{\partial \Phi_j}{\partial \sigma_j}|_{\sigma_j=0}=r_j(U_b)$.
Set
$$
\Phi(\sigma_5,\sigma_4,\sigma_3,\sigma_2,\sigma_1;U_b)=\Phi_5(\sigma_5,\Phi_4(\sigma_4,\Phi_3(\sigma_3,\Phi_2(\sigma_2,
\Phi_1(\sigma_1;U_b))))).
$$
We denote $\Psi_j(\sigma_j,W(U_b))=W(\Phi_j(\sigma_j;U_b))$ and
\begin{eqnarray*}
\Psi(\sigma_5,\sigma_4,\sigma_3,\sigma_2,\sigma_1;W(U_b))
&=&\Psi_5(\sigma_5,\Psi_4(\sigma_4,\Psi_3(\sigma_3,\Psi_2(\sigma_2,\Psi_1(\sigma_1
;W(U_b))))))\\
&=&W(\Phi(\sigma_5,\sigma_4,\sigma_3,\sigma_2,\sigma_1;U_b)).
\end{eqnarray*}

\par Finally, we denote
\begin{equation}
\boldsymbol{\sigma}=(\sigma_1,\sigma_2,\sigma_3,\sigma_4,\sigma_5),
\end{equation}
 and
 \begin{equation}
 \Psi(\boldsymbol{\sigma},W(U_b))=\Psi_5(\sigma_5,\Psi_4(\sigma_4,\Psi_3(\sigma_3,\Psi_2(\sigma_2,\Psi_1(\sigma_1,W(U_b)))))).
\end{equation}

\section{The Glimm Fractional-Step Scheme}
We employ a fractional-step scheme for the inhomogeneous system \eqref{2}
as described in \cite{chen1} based on the Glimm scheme.
As before, we regard the $x$-direction as the time-like direction.

\par Choose mesh lengths $h>0$ and $l>0$ in the $x$-direction and $y$-direction, respectively,
such that the Courant-Friedrichs-Levy condition holds:
\begin{equation}
\Lambda=\max_{1\le j \le 5}|\lambda_j(U)|\le \frac{l}{2h}.
\label{c2}
\end{equation}

\par Partition $\mathbb{R}^+$ by the sequence $x_k=kh,k\in \mathbb{Z}^+$,
and partition $\mathbb{R}$ into cells with the $j$th cell centered at
\begin{equation*}
 y_j=jl, \qquad j=0,\pm 1,\pm 2, \cdots.
\end{equation*}
We begin with approximating the initial data $U_0(y)$ by a function $U^h(0,y)$,
which is constant for $y$ in the interval $[y_{j-1}, y_{j+1}]$ for $j$ even
and converges to $U_0(y)$ both pointwise a.e. and in $L^1$ on any bounded
interval as $h\to 0$.
Choose a random sequence $\theta_k, k=0,1,2,\cdots$, in the interval $(-1,1)$ with
the uniform probability distribution.

\medskip
\par We then construct the approximate solution $W(U^h(x,y))$ as follows:

\par Assume that $W(U^h(x,y))$ is defined for $x<kh$.
Then we construct the approximate solution $W(U^h(x,y))$ in the
strip $[kh,(k+1)h)\times (-\infty, \infty)$ as follows:

\medskip
\par \textit{Step 1  (Random step):} Define
\begin{eqnarray*}
&&W(U_j^k)=W(U^h(kh-,(j+\theta_k)l)),\\[1.5mm]
&&W(U^h(kh+0,y))\equiv W(U_j^k), \qquad (j-1)l\le y <(j+1)l,
\end{eqnarray*}
where $j+k$ is even, and $\chi _k $ is the $k$th element of the random
sequence $(\chi_1, \cdots,\chi_k,\cdots)$.

\medskip
\par \textit{Step 2  (Solving the Riemann problem):}
In the strip $[kh,(k+1)h)\times (-\infty,\infty)$, we solve the following Riemann
problem in each domain $(kh,(k+1)h)\times ((j-1)l,(j+1)l)$:
\begin{eqnarray*}
\left\{
\begin{array}{ll}
W(U)_{\tau}+H(U)_y=0,\\[2mm]
 W(U)|_{\tau=0}=
\begin{cases} W(U_{j-1}^k) & y< jl,\\[1.5mm]
W(U_{j+1}^k) & y> jl,
\end{cases}
\end{array} \right.
\end{eqnarray*}
where $j+k$ is odd and $\tau =x-kh$. The resulting solution is denoted as $W(U_0^h(x,y))$.
\begin{center}
\setlength{\unitlength}{1mm}
\begin{picture}(48,60)(0,-60)
\linethickness{1pt}
\multiput(0,0)(0,-3){16}{\line(0,-2){2}}
\multiput(24,0)(0,-3){16}{\line(0,-2){2}}
\put(0,-24){\line(5,4){20}}
\thicklines
\put(0,-24){\line(4,1){20}}
\thinlines
\put(0,-24){\line(3,-1){20}}
\put(8,-4){$W(U_{j+1}^k)$}
\put(8,-40){$W(U_{j-1}^k)$}
\put(12,-10){$\sigma_5$}
\put(14,-17){$\sigma_{2(3,4)}$}
\put(15,-28){$\sigma_1$}
\put(-4,-60){Fig. 2. Riemann problem}
\put(-4,-52){$x=kh$}
\put(20,-52){$(k+1)h$}
\put(-14,-26){$(kh,jl)$}
\end{picture}
\end{center}

\par \textit{Step 3  (Reacting step):} Define
\begin{equation*}
W(U^h(x,y))=W(U_0^h(x,y))+G(U_0^h(x,y))(x-kh),
\end{equation*}
where $G(U)=(0,0,0,q\rho Z\phi(T),-\rho Z\phi(T))$ as before and $kh\le x<(k+1)h$.

\par Therefore, we can construct the approximate solution $W(U^h(x,y))$
in the strip $[kh,(k+1)h)\times (-\infty, \infty)$ as long as the Riemann problems
in Step 2 are solvable.

\section{$BV$--Stability}

\par In this section, we estimate the approximate solutions $W(U^h(x,y))$ in the total
variation norm and prove that the total variation of
the approximate solutions $W(U^h(x,y))$ in $y$, for any fixed $x$,
is uniformly bounded with respect to the mesh
length $h$. We measure the total variation of approximate solutions by using the sum of the
absolute values of the strengths of waves in the solution of each Riemann problem
in Step 2 as in Section 3.

\par We define a weighted $l_1$--norm
\begin{equation}\label{4.1a}
\|v\|_1=|v_1|+|v_2|+|v_3|+M|v_4|+|v_5| \qquad \text{for a vector $v=(v_1,v_2,v_3,v_4,v_5) \in \mathbb{R}^5$},
\end{equation}
where $M>0$ is a constant to be determined later.

\par We define another norm
\begin{equation}
\|v\|=|v_1|+|v_2|+|v_3|+|v_4| \qquad \text{for a vector $v=(v_1,v_2,v_3,v_4)\in \mathbb{R}^4$}.
\end{equation}

\par Let $U_g\equiv (u,v,p,\rho)$ and $W_g\equiv (\rho u,\rho u^2+p,\rho uv, \rho u(\bar{h}+\frac{u^2+v^2}{2}))$
denote the first four components of $U$ and $W$, respectively.

\subsection{\textbf{Interaction estimates on the non-reacting step}}

The interaction estimate for \eqref{3} is similar to the argument for Proposition 3.1 in \cite{chen2}.

\begin{center}
\setlength{\unitlength}{1mm}
\begin{picture}(50,74)(-2,-4)
\linethickness{1pt}
\multiput(0,0)(0,3){24}{\line(0,2){2}}
\multiput(24,0)(0,3){24}{\line(0,2){2}}
\multiput(48,0)(0,3){24}{\line(0,2){2}}
\put(0,54){\line(3,2){20}}
\put(0,54){\line(3,-1){20}}
\put(0,54){\line(5,2){20}}
\put(0,18){\line(1,1){20}}
\put(0,18){\line(5,3){20}}
\put(0,18){\line(4,-1){20}}
\put(24,40){\line(5,-1){20}}
\put(24,40){\line(3,2){20}}
\put(24,40){\line(5,1){20}}
\put(2,70){$U_a$}
\put(2,30){$U_m$}
\put(2,5){$U_b$}
\put(36,15){$U_b$}
\put(36,60){$U_a$}
\put(12,66){$\beta_5$}
\put(12,56){$\beta_{2(3,4)}$}
\put(12,46){$\beta_1$}
\put(12,36){$\alpha_5$}
\put(12,23){$\alpha_{2(3,4)}$}
\put(12,11){$\alpha_1$}
\put(40,54){$\gamma_5$}
\put(38,46){$\gamma_{2(3,4)}$}
\put(40,38){$\gamma_1$}
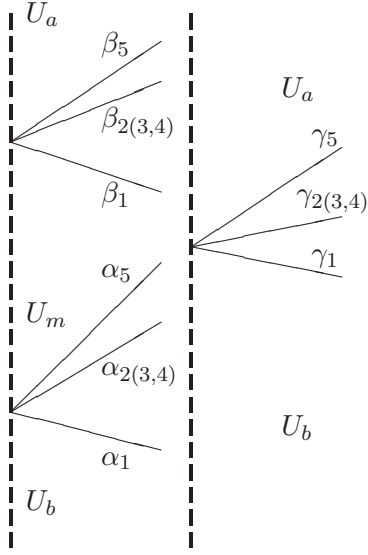
\put(-1,-6){Fig. 7. Weak wave interaction}
\end{picture}
\end{center}

\medskip
\begin{lemma}\label{4.1b}
Suppose that $U_b, U_m$, and $U_a$ are three states
in a small neighborhood $O_{\varepsilon}(U_+)$ with
\begin{equation*}
\{U_b, U_m\}=(\alpha_1,\alpha_2,\alpha_3,\alpha_4,\alpha_5),\{U_m,U_a\}=(\beta_1,\beta_2,\beta_3,\beta_4,\beta_5),
\{U_b,U_a\}=(\gamma_1,\gamma_2,\gamma_3,\gamma_4,\gamma_5).
\end{equation*}
Then
\begin{equation*}\label{33}
\gamma_i=\alpha_i+\beta_i+O(1)\Delta(\alpha,\beta),
\end{equation*}
where
$\Delta(\alpha,\beta)=|\alpha_5|(|\beta_1|+|\beta_2|+|\beta_3|+|\beta_4|)+|\beta_1|(|\alpha_2|+|\alpha_3|+|\alpha_4|)
+\sum_{j=1,5}\Delta_j(\alpha,\beta)$ with
\begin{equation*}
\Delta_j(\alpha,\beta)=\left\{
        \begin{array}{ll}
        0,&\quad \mbox{$\alpha_j \ge 0$ and $\beta_j \ge 0$},\\
        |\alpha_j||\beta_j|,&\quad \mbox{otherwise.}
        \end{array}
                        \right.
\end{equation*}
\end{lemma}

\subsection{\textbf{Estimates on the reacting step}}

\par For convenience, we use $\tilde{U}$ to denote the value of $U$ before reaction,
while $U$  after reaction. That is,
\begin{equation*}
W(U(x,y))=W(\tilde{U}(x,y))+G(\tilde{U}(x,y))\tau,
\end{equation*}
where $\tau = x-kh$ and $kh\le x <(k+1)h$.

\begin{lemma} Let
\begin{equation}
\boldsymbol{\sigma}=(\sigma_1,\sigma_2,\sigma_3,\sigma_4,\sigma_5)=B(W(\tilde{U}_b),W(\tilde{U}_a))
\end{equation}
be the vector of signed wave strengths in the solution of the Riemann problem
with Riemann data $(W(\tilde{U}_b),W(\tilde{U}_a))$.
Let
\begin{equation}
\Gamma(W(\tilde{U}_b),\boldsymbol{\sigma},h)=B(W(U_b),W(U_a)),
\label{a4}
\end{equation}
where $W(U_b)=W(\tilde{U}_b)+G(\tilde{U}_b)h$ and $W(U_a)=W(\tilde{U}_a)+G(\tilde{U}_a)h$.
Then
\begin{equation}
\Gamma(W(\tilde{U}_b),\boldsymbol{\sigma},h)=\boldsymbol{\sigma}+O(||\boldsymbol{\sigma}||_1)h.
\end{equation}
\end{lemma}

\par Lemma 4.2 implies that the increasing of the total variation of the fractional-step
approximate solutions
is at no more than an exponential rate.

\medskip
\begin{lemma}
$\|\Gamma(W(\tilde{U}_b),\boldsymbol{\sigma},h)
-\big(\boldsymbol{\sigma}+\frac{\partial \Gamma}{\partial h}(W(\tilde{U}_b),\boldsymbol{\sigma},0)h\big)\|_1
\le C\|\boldsymbol{\sigma}\|_1\frac{h^2}{2}.
$
\end{lemma}

\par Lemma 4.3 shows that we can estimate the increase in the total variation for the reacting step
by calculating the first derivatives of the solution operator for the Riemann problem.

\par The proof of Lemmas 4.2--4.3 can be found in \cite{chen1}.

\medskip
In particular, for \eqref{d1}--\eqref{d2}, we need to analyze the reacting step, which takes the form:
\begin{equation*}
W(U^h(x,y))=W(U_0^h(x,y))+G(U_0^h(x,y))\tau,
\end{equation*}
where all the quantities $\rho_0^h, \rho^h$, {\it etc.} are evaluated at $(kh+\tau,y)$ and $\tau =x-kh$.
More precisely, it takes the form:
\begin{equation}
\begin{split}
\rho^hu^h&=\rho_0^hu_0^h,\\[1.5mm]
\rho^h(u^h)^2+p^h&=\rho_0^h(u_0^h)^2+p_0^h,\\[1.5mm]
\rho^hu^hv^h&=\rho_0^hu_0^hv_0^h,\\[1.5mm]
(\rho^hE^h+p^h)u^h&=(\rho_0^hE_0^h+p_0^h)u_0^h+q\rho_0^hz_0^h\phi(T_0^h)\tau,\\[1.5mm]
\rho^hu^hZ^h&=\rho_0^hu_0^hZ_0^h-\rho_0^hZ_0^h\phi(T_0^h)\tau.
\end{split} \label{diff}
\end{equation}
We need to estimate the change in $(u,v,p,\rho,z, T)$ due to the reaction step.

\par First, we have
\begin{equation*}
(T^h-T_0^h)(kh+\tau,y)=\frac{\partial T}{\partial w_4}q\rho_0^hZ_0^h\phi(T_0^h)\tau
=\frac{(\gamma -1)((u_0^h)^2-RT_0^h)}{R\rho_0^hu_0^h((u_0^h)^2-\gamma RT_0^h)}q\rho_0^hZ_0^h\phi(T_0^h)\tau.
\end{equation*}
Since $u^2>c^2=\frac{\gamma p}{\rho}=\gamma RT$, then $T^h(x,y)\ge T_0^h(x,y)$,
which shows that the temperature $T$ does not decrease due to the reaction.

\par Second, from the fifth equation: $Z^h-Z_0^h=-\frac {Z_0^h\phi(T_0^h)\tau}{u_0^h}$.
Since $\phi(T)$ is assumed to be Lipschitz continuous, nonnegative, and increasing,
there exists a constant $\Phi_1>0$ such that
\begin{equation}
Z^h-Z_0^h\le -Z_0^h\Phi_1 \tau.
\end{equation}
Then we conclude
\begin{equation}
Z^h\le Z_0^h(1-\Phi_1 \tau)\le Z_0^h e^{-\Phi_1 \tau}, \qquad 0\le \tau<h.
\label{a1}
\end{equation}
According to the scheme and using the induction,
we can actually obtain
\begin{equation}
Z_0^h\le \|Z_0\|_{\infty}e^{-\Phi_1 kh}, \qquad kh \le x<(k+1)h.
\end{equation}

\par Third, from the first three equations,
we know that $u^h=\frac{\rho_0^hu_0^h}{\rho^h}$, $v^h=v_0^h$,
and $p^h=p_0^h+ \rho_0^h(u_0^h)^2-\frac{\rho_0^h(u_0^h)^2}{\rho^h}$.
Substitution of these into the fourth equation, we have
\begin{equation}
\frac{\gamma+1}{2}(\rho_0^h u_0^h)^2(\frac{1}{\rho^h})^2
-\gamma\big(\rho_0^h(u_0^h)^2+p_0^h\big)\frac{1}{\rho^h}
+\big(\frac{\gamma-1}{2}(u_0^h)^2+\gamma \frac{p_0^h}{\rho_0^h}+Z_0^hO(h)\big)=0.
\end{equation}
Therefore, we obtain
\begin{equation}
\frac{1}{\rho^h}=\frac{\gamma\big(\rho_0^h(u_0^h)^2+p_0^h\big)
+\sqrt{\big(\rho_0^h(u_0^h)^2-\gamma p_0^h\big)^2+(\rho_0^h u_0^h)^2Z_0^hO(h)}}{(\gamma+1)(\rho_0^hu_0^h)^2}.
\end{equation}
Using the Taylor expansion, we know
\begin{equation}
\frac{1}{\rho_0^h}=\frac{1}{\rho^h}+Z_0^hO(h).
\end{equation}
That is,
\begin{equation}
\rho^h-\rho_0^h=\|Z_0\|_{\infty}e^{-\Phi_1 kh}O(h).
\end{equation}
\par Similar calculations also apply to $u$ and $p$.
Therefore, we have

\begin{lemma}
There are positive constants $C_0$ and $\Phi_1$ such that
\begin{equation}
\begin{split}
T^h \ge T_0^h &\ge C_0>0,\\
u^h-u_0^h&=\|Z_0\|_{\infty}e^{-\Phi_1 kh} O(h),\\
v^h-v_0^h&=0,\\
p^h-p_0^h&=\|Z_0\|_{\infty}O(h)e^{-\Phi_1 kh},\\
\rho^h-\rho_0^h&=\|Z_0\|_{\infty}O(h)e^{-\Phi_1 kh},\\
Z^h &\le Z_0^h e^{-\Phi_1 \tau}, \qquad 0\le \tau<h.
\end{split}
\end{equation}
Furthermore, \begin{equation}
Z_0^h\le \|Z_0\|_{\infty}e^{-\Phi_1 kh}O(h), \qquad kh \le x<(k+1)h.
\end{equation}
All the quantities are evaluated at $(kh+\tau,y)$ with $\tau =x-kh$.
\label{q1}
\end{lemma}

\subsection{\textbf{Glimm functional for the fractional-step scheme}}

\par Following Glimm's method \cite{glimm}, we define a functional on the restriction
of the approximate solution $W(U^h)$ to certain mesh curves $J$.
We define a mesh point to be a point $(x,y)=(kh,(j+\theta_k)l)$,
where $k\in \mathbb{N}$ and $j\in \mathbb{Z} $ such that $j+k$ is even.
A mesh curve $J$ is a piecewise linear curve in the $(x,y)$--plane,
which successively connects the mesh points $(kh,(j+\theta_k)l)$ to the mesh points
$((k\pm 1)h,(j+1+\theta_{k\pm 1})l)$.
We define a partial order on the set of mesh curves by stating that larger
curves lie toward larger $x$.
We call $J_2$ an immediate successor of $J_1$ if $J_2$ connects the same mesh points
as  $J_1$, except for one mesh point, and if $J_2>J_1$.
Let $J_k$ be the unique mesh curve which connects the mesh points on $x=kh$
to the mesh points on $x=(k+1)h$. Note that $J_k$ crosses all the waves
in the Riemann solutions of $W(U_0^h(x,y))$ in the strip $kh\le x<(k+1)h$.

\par We now define a functional $F$ on the set of mesh curves.
For any mesh curve $J$, we define
\begin{equation}
L_i(J)=\sum\{|\alpha|:\text{$\alpha$ is the $i$th wave crossing $J$}\}
\qquad\mbox{for $1\le i \le 5$}.
\end{equation}

\par Next, we define
\begin{equation}
L(J)=\sum_{1\le i\le 5,i\ne 4}L_i(J)+M L_4(J),
\end{equation}
and
\begin{equation}
Q(J)=\sum \{|\alpha| |\beta|: \text{both $\alpha$ and  $\beta$ cross $J$ and approach each other}\},
\end{equation}
where $M>0$ is a constant to be determined as in \eqref{4.1a}.

By standard procedure as in \cite{smoller1}, when $\mathrm{TV}(U_0(\cdot))$ is small enough,
we can choose a positive constant $K_0$ sufficiently large such that
the Glimm functional
\begin{equation}
F(J)=L(J)+ K_0Q(J)
\end{equation}
is non-increasing in the non-reacting step.

\subsection{\textbf{BV-stability of the reaction step}}

\par We now prove the BV-stability of the approximate solutions during the reaction step.
Our total variation bounds imply bounds on the length of $W(U^h(J))$,
but we must also deal with the ``drift'' of the solution due to the reaction term $G(U)$.

\par In order to discuss the effect of the exothermic reaction on the functionals $L$
and $Q$, it is convenient to identify a new ``mesh curve" $\tilde J$, which, as a curve,
is the same as a given mesh curve $J$, but upon which the value of $W(U)$ differs from
the value of $W(U)$ on $J$ by a single reaction step along all of $J$.
We take $\tilde J$ to represent the values before the reaction and $J$ to represent
the values after the reaction step.

\begin{lemma}
There is a positive constant C such that
\begin{gather}
L(J_k)\le L(\tilde{J_k})+ Cqh\|w_{5,0}\|_{\infty}e^{-\Phi_1 kh}L(\tilde{J_k}),\\[2mm]
Q(J_k)\le Q(\tilde{J_k})+ Cqh \|w_{5,0}\|_{\infty}e^{-\Phi_1kh}L(\tilde{J_k})^2.
\end{gather}
\end{lemma}

\begin{pf}
For simplicity of presentation, we denote $\textbf{c}=(0,0,0,1,-\frac{1}{q})^\top, \boldsymbol{\tilde{\sigma}}_i=(\tilde{\sigma}_{1i},\tilde{\sigma}_{2i},\tilde{\sigma}_{3i},\tilde{\sigma}_{4i},
\tilde{\sigma}_{5i})$, $\textbf{c}_g=(0,0,0,1)^\top$,
$W(\tilde{U}_{i+1})=\Psi(\boldsymbol{\tilde{\sigma}}_i,W(\tilde{U_i}))$,
$B=(B_1,B_2,B_3,B_4,B_5)^\top$, and $B_g=(B_1,B_2,B_3,B_4)^\top$.

\par Let
\begin{equation}
\Gamma(W(\tilde{U_i}),\boldsymbol{\tilde{\sigma}_i},h)=B(W(U_i),W(U_{i+1}))
\end{equation}
as before, where $W(U_i)=W(\tilde{U_i})+G(\tilde{U_i})h$
and $W(U_{i+1})=W(\tilde{U}_{i+1})+G(\tilde{U}_{i+1})h$. Then we have
\begin{equation*}
\begin{split}
\frac{\partial \Gamma}{\partial h}(W(\tilde{U_i}),\boldsymbol{\tilde{\sigma}}_i,0)
&=\tilde{\rho}_{i}\tilde{Z}_{i}\phi(\tilde{T}_{i}){\partial}_1Bq\textbf{c}
 +\tilde{\rho}_{i+1}\tilde{Z}_{i+1}\phi(\tilde{T}_{i+1}){\partial}_2Bq\textbf{c}\\
&=\tilde{w}_{5,i}\frac{\phi(\tilde{T}_i)}{\tilde{u}_i}{\partial}_1Bq\textbf{c}
+\tilde{w}_{5,i+1}\frac{\phi(\tilde{T}_{i+1})}{\tilde{u}_{i+1}}{\partial}_2Bq\textbf{c}\\
&=\tilde{w}_{5,i}({\partial_1} B \frac{\phi(\tilde{T}_i)}{\tilde{u}_i}+{\partial}_2 B\frac{\phi(\tilde{T}_{i+1})}{\tilde{u}_{i+1}})q\textbf{c}+(\tilde{w}_{5,i+1}-\tilde{w}_{5,i}){\partial}_2 B\frac{\phi(\tilde{T}_{i+1})}{\tilde{u}_{i+1}}q\textbf{c},
\end{split}
\end{equation*}
where $\tilde{w}_{5,i}=\tilde{\rho}_{i}\tilde{u}_{i}\tilde{Z}_{i}$
and $\tilde{w}_{5,i+1}=\tilde{\rho}_{i+1}\tilde{u}_{i+1}\tilde{Z}_{i+1}$.
Since $Z$ is decoupled from  $(u,v,p,\rho)^\top$ in the solution of the non-reacting Riemann problem,
this means that ${\partial}_1B$ and ${\partial}_2B$
are the block $5\times 5$ matrices with the upper left $4\times 4$ block relating to non-reacting gas dynamics.
The remaining $1\times 1$ block contains the derivative of wave strength of $Z$-contact with respect
to $\rho uZ$---the value of this derivative is $-1$ for $\frac{\partial B_5}{\partial w_{5,i}}$
and $1$ for $\frac{\partial B_5}{\partial w_{5,i+1}}$,
since $B_5=\rho_{i+1} u_{i+1} Z_{i+1}-\rho_{i} u_{i} Z_{i}=w_{5,i+1}-w_{5,i}$.
Then we have
\begin{equation}
{\partial}_1 Bq\textbf{c}=\left(\begin{array}{ccc}
{\partial}_{1W_g} B_g&0\\
0&-1
\end{array}\right)
\left(\begin{array}{ccc}
q\textbf{c}_g\\
-1
\end{array}\right)
=
\left(\begin{array}{ccc}
{\partial}_{1W_g} B_gq\textbf{c}_g\\
0
\end{array}\right)+(0,0,0,0,1)^\top,
\end{equation}
\begin{equation}
{\partial}_2 Bq\textbf{c}=\left(\begin{array}{ccc}
{\partial}_{2W_g} B_g&0\\
0&1
\end{array}\right)
\left(\begin{array}{ccc}
q\textbf{c}_g\\
-1
\end{array}\right)
=
\left(\begin{array}{ccc}
{\partial}_{2W_g} B_gq\textbf{c}_g\\
0
\end{array}\right)-(0,0,0,0,1)^\top,
\end{equation}
where $W_g=(w_1, \cdots, w_4)$. Then
\begin{equation}
\begin{split}
\frac{\partial \Gamma}{\partial h}(W(\tilde{U_i}),\boldsymbol{\tilde{\sigma}}_i,0)
=&\tilde{w}_{5,i}\Bigg[\frac{\phi(\tilde{T}_{i})}{\tilde{u}_{i}}\left(\begin{array}{ccc}
{\partial}_{1W_g} B_gq\textbf{c}_g\\
1
\end{array}\right)
+ \frac{\phi(\tilde{T}_{i+1})}{\tilde{u}_{i+1}}\left(\begin{array}{ccc}
{\partial}_{2W_g} B_gq\textbf{c}_g\\
-1
\end{array}\right)\Bigg]\\&+(\tilde{w}_{5,i+1}-\tilde{w}_{5,i})\frac{\phi(\tilde{T}_{i+1})}{\tilde{u}_{i+1}}\left(\begin{array}{ccc}
{\partial}_{2W_g} B_gq\textbf{c}_g\\
-1
\end{array}\right).
\label{b1}
\end{split}
\end{equation}

Thus, the first four components of \eqref{b1} have the form:
\begin{equation}
\begin{split}
\frac{\partial \Gamma_g}{\partial h}(W(\tilde{U}_i),\boldsymbol{\tilde{\sigma}}_i,0)
=&\tilde{w}_{5,i}\Big[{\partial}_{1W_g} B_gq\textbf{c}_g\frac{\phi(\tilde{T}_{i})}{\tilde{u}_{i}}
 + {\partial}_{2W_g} B_gq\textbf{c}_g\frac{\phi(\tilde{T}_{i+1})}{\tilde{u}_{i+1}}\Big]\\
 &+(\tilde{w}_{5,i+1}-\tilde{w}_{5,i}){\partial}_{2W_g} B_gq\textbf{c}_g\frac{\phi(\tilde{T}_{i+1})}{\tilde{u}_{i+1}}\\
=& \tilde{w}_{5,i}A(W(\tilde{U}_i),\boldsymbol{\tilde{\sigma}}_i)+(\tilde{w}_{5,i+1}
-\tilde{w}_{5,i}){\partial}_{2W_g} B_gq\textbf{c}_g \frac{\phi(\tilde{T}_{i+1})}{\tilde{u}_{i+1}},
\label{b2}
\end{split}
\end{equation}
where
\begin{equation}
A(W(\tilde{U}_i),\boldsymbol{\tilde{\sigma}}_i)={\partial}_{1W_g} B_gq\textbf{c}_g\frac{\phi(\tilde{T}_{i})}{\tilde{u}_{i}}
+ {\partial}_{2W_g} B_gq\textbf{c}_g\frac{\phi(\tilde{T}_{i+1})}{\tilde{u}_{i+1}}.
\end{equation}
It is easy to see that, if $\boldsymbol{\tilde{\sigma}}_{g,i}=(\tilde{\sigma}_{1,i},\tilde{\sigma}_{2,i},\tilde{\sigma}_{3,i},\tilde{\sigma}_{4,i})=\boldsymbol{0}$,
then $\tilde{W}_{g,i+1}=\tilde{W}_{g,i}:=W_g(\tilde{U}_i)$ and, in particular, $\tilde{T}_{i+1}=\tilde{T}_i$.
Since $B_g(\tilde{W}_{g,i},\tilde{W}_{g,i})$ is the vector of wave strengths for a Riemann problem with equal states,
\begin{equation}
{\partial}_{1W_g} B_g|_{\boldsymbol{\tilde{\sigma}}_{g,i}=0}={\partial}_{2W_g} B_g|_{\boldsymbol{\tilde{\sigma}}_{g,i}=0}=0.
\end{equation}
Therefore, there exists some positive constant $C$ such that the first term in \eqref{b2} can be estimated by
\begin{equation}
 \|\tilde{w}_{5,i}A(W(\tilde{U}_i),\boldsymbol{\tilde{\sigma}}_i)\|
 \le C\tilde{w}_{5,i} \|\boldsymbol{\tilde{\sigma}}_{g,i}\|q.
\end{equation}

\par We next examine the last term of \eqref{b2}, which has the form
\begin{equation}
(\tilde{w}_{5,i+1}-\tilde{w}_{5,i}){\partial}_{2W_g} B_gq\textbf{c}_g \frac{\phi(\tilde{T}_{i+1})}{\tilde{u}_{i+1}}.
\label{b3}
\end{equation}
The fifth component of \eqref{b1} is the equation for the strength of the $Z$-wave. This equation is
\begin{equation*}
\frac{\partial}{\partial h}(\tilde{w}_{5,i+1}-\tilde{w}_{5,i})
=\tilde{w}_{5,i}(\frac{\phi(\tilde{T}_i)}{\tilde{u}_i}-\frac{\phi(\tilde{T}_{i+1})}{\tilde{u}_{i+1}})
-(\tilde{w}_{5,i+1}-\tilde{w}_{5,i})\frac{\phi(\tilde{T}_{i+1})}{\tilde{u}_{i+1}},
\end{equation*}
so that
\begin{equation*}
\frac{\partial}{\partial h}|\tilde{w}_{5,i+1}-\tilde{w}_{5,i}|
\le \tilde{w}_{5,i}|\frac{\phi(\tilde{T}_i)}{\tilde{u}_i}
-\frac{\phi(\tilde{T}_{i+1})}{\tilde{u}_{i+1}}|-|\tilde{w}_{5,i+1}-\tilde{w}_{5,i}|\frac{\phi(\tilde{T}_{i+1})}{\tilde{u}_{i+1}}.
\end{equation*}

\par Thus, the reaction step produces possible increases in the total variation, which are bounded by
\[
C\tilde{w}_{5,i}\|\boldsymbol{\tilde{\sigma}}_{g,i}\|qh+|\tilde{w}_{5,i+1}
-\tilde{w}_{5,i}|\frac{\phi(\tilde{T}_{i+1})}{\tilde{u}_{i+1}}\|{\partial}_{2W_g} B_g\textbf{c}_g\|qh.
\]
The reaction step also produces a decrease in total variation for the $w_5=\rho uZ$
component---the fifth component of
$\frac{\partial \Gamma}{\partial h}(W(\tilde{U}_{i}),\boldsymbol{\tilde{\sigma}}_i,0)$---in the amount $|\tilde{w}_{5,i+1}-\tilde{w}_{5,i}|\frac{\phi(\tilde{T}_{i+1})}{\tilde{u}_{i+1}}$.
We now use the decrease in the $w_5$--component proportional to $|\tilde{w}_{5,i+1}-\tilde{w}_{5,i}|$.
Since ${\partial}_{2W_g} B_g$ is Lipschitz continuous, there exists a upper bound $M$ for $\|{\partial}_{2g} B_g\|q$.
Thus, the effect of term \eqref{b2} on $(\rho u,\rho u^2+p,\rho uv, \rho u (\bar{h}+\frac{u^2+v^2}{2}))$ of
$\frac{\partial \Gamma}{\partial h}(W(\tilde{U}_{i}),\boldsymbol{\tilde{\sigma}}_i,0)$ is bounded by $M|\tilde{w}_{5,i+1}-\tilde{w}_{5,i}|\frac{\phi(\tilde{T}_{i+1})}{\tilde{u}_{i+1}}h$,
and this increase is offset by a decrease in the term $M|\tilde{w}_{5,i+1}-\tilde{w}_{5,i}|$.

\par Thus, the change in $L$ is estimated as follows:
\begin{equation}
\begin{split}
L(J_k)- L(\tilde{J_k})&=\sum_{1\le j \le 5,j\ne 4}\big(L_j(J_k)- L_j(\tilde{J_k})\big)+M\big(L_4(J_k)- L_4(\tilde{J_k})\big)
\\&= \sum_{1\le j \le 5,j\ne 4}\sum_{-\infty <i<\infty}(|\sigma_{j,i}|-|\tilde{\sigma_{j,i}}|)+M \sum_{-\infty <i<\infty}(|\sigma_{4,i}|-|\tilde{\sigma_{4,i}}|)
\\& \le \sum_{-\infty <i<\infty}\|\frac{\partial \Gamma}{\partial h}(W(\tilde{U}_i),\boldsymbol{\tilde{\sigma}}_i,0)\|_1h
\\& \le \sum_{-\infty <i<\infty} \Big(C q\tilde{w}_{5,i}\|\boldsymbol{\tilde{\sigma}}_{g,i}\|h
 +|\tilde{w}_{5,i+1}-\tilde{w}_{5,i}|\frac{\phi(\tilde{T}_{i+1})}{\tilde{u}_{i+1}}q\|{\partial}_{2g} B_g\textbf{c}_g\|h
\\&
\qquad \qquad\qquad
+M \big(\tilde{w}_{5,i}|\frac{\phi(\tilde{T}_i)}{\tilde{u}_i}-\frac{\phi(\tilde{T}_{i+1})}{\tilde{u}_{i+1}}|h
-|\tilde{w}_{5,i+1}-\tilde{w}_{5,i}|\frac{\phi(\tilde{T}_{i+1})}{\tilde{u}_{i+1}}h \big)\Big)
\\&\le Cqh \|\tilde{w_5}\|_{\infty}L(\tilde{J_k})
\\&\le Cqh\|w_{5,0}\|_{\infty}e^{-\Phi_1kh}L(\tilde{J_k}),
\end{split}
\end{equation}
where we have chosen $M>0$ large enough to make the third inequality hold, and the last inequality comes from Lemma \ref{q1}.

\par Consequently, we have
\begin{equation}
\begin{split}
 Q(J_k)-Q(\tilde{J_k})&=\sum_{App}\big(|\alpha_i||\beta_j|-|\tilde{\alpha_i}||\tilde{\beta_j}|\big)
 \\&= \sum_{App}\big(|\alpha_i|(|\beta_j|- |\tilde{\beta_j}|)+|\tilde{\beta_j}|( |\alpha_i|- |\tilde{\alpha_i}|)\big)
 \\&\le C\big(L(J_k)- L(\tilde{J_k})\big)L(\tilde{J_k}).
 \end{split}
\end{equation}
Therefore,
 \begin{equation}
 \begin{split}
Q(J_k)\le Q(\tilde{J_k})+Cqh\|w_{5,0}\|_{\infty}e^{-\Phi_1kh}L(\tilde{J_k})^2
\le Q(\tilde{J_k})+Cqh\|w_{5,0}\|_{\infty}e^{-\Phi_1kh}F(\tilde{J_k})^2.
\end{split}
\end{equation}
The proof is completed.
\end{pf}

\par Since
\begin{equation}
F(J)=L(J)+ K_0Q(J),
\end{equation}
we have actually proved the following lemma.

\begin{lemma}
Let $J_k$ be a mesh curve between $x=kh$ and $x=(k+1)h$. Then
\begin{equation}
F(J_k)\le F(\tilde{J_k})\big(1+Cqh\|w_{5,0}\|_{\infty}e^{-\Phi_1 kh}(1+ F(\tilde{J_k}))\big),
\label{11}
\end{equation}
where C is a constant independent of the mesh lengths $l$ and $h$.
\end{lemma}

\par We need to obtain a uniform bound on $F$.
First of all, we suppose that such a bound exists, namely, $F(\tilde{J}_k)\le A$
for some positive constant.  Then, by (\ref{11}), we have
\[
F(J_k)\le F(\tilde{J_k})\big(1+Cq\|w_{5,0}\|_{\infty}e^{-\Phi_1 kh}(1+A)h\big).
\]
Since $F$ is non-increasing in the non-reacting step, $F(\tilde{J}_k)\le F(J_{k-1})$.
Then we have
\[
F(J_k)\le F(\tilde{J}_0)\prod_{j=0}^{k}\big(1+Cq\|w_{5,0}\|_{\infty}d^{j}(1+A)h\big),
\]
where $d=e^{-\Phi_1h}$.
Using the inequality $\mathrm{ln}(1+x)\le x$ for $x\ge 0$,
\begin{equation}
\begin{split}
\mathrm{ln}\Big(\frac{F(J_k)}{F(\tilde{J}_0)}\Big)
&\le \sum_{j=0}^k \mathrm{ln}\big(1+Cq\|w_{5,0}\|_{\infty}d^j(1+A)h\big)
\\&\le \sum_{j=0}^k Cqh \|w_{5,0}\|_{\infty}d^j(1+A)
\\&\le Cqh \|w_{5,0}\|_{\infty}(1+A)\frac{1}{1-d}.
\end{split}
\end{equation}
Thus we obtain
\begin{equation}
F(J_k)\le F(\tilde{J}_0)\mathrm{exp}\Big(\frac{Cqh\|w_{5,0}\|_{\infty}(1+A)}{1-e^{-\Phi_1 h}}\Big).
\label{15}
\end{equation}
The function $f(h)=\frac{h}{1-e^{-\Phi_1 h}}$ is increasing for $h>0$ and tends to $\frac{1}{\Phi_1}$
as $h\to 0$. Thus, for $h$ sufficiently small, we obtain
\begin{equation}
F(J_k)\le F(\tilde{J}_0)\mathrm{exp}\Big(\frac{C_1q\|w_{5,0}\|_{\infty}(1+A)}{\Phi_1}\Big),
\label{12}
\end{equation}
where $C_1=2C$.
Estimate (\ref{12}) is valid as long as $F(\tilde{J}_k)\le A$.
Since $F(\tilde{J}_k)\le F(J_{k-1})$, the condition required for this result is that
\begin{equation}
F(\tilde{J}_0)\le \mathrm{exp}\Big(-\frac{C_1q\|w_{5,0}\|_{\infty}(1+A)}{\Phi_1}\Big)A=:g(A).\label{13}
\end{equation}
The value of A which maximizes $g(A)$ is $A=\frac{\Phi_1}{C_1q\|w_{5,0}\|_{\infty}}$.
Thus, our least-restrictive condition on $F(\tilde{J}_0)$ is
\begin{equation}
F(\tilde{J}_0)
\le \mathrm{exp}\Big(-1-\frac{C_1q\|w_{5,0}\|_{\infty}}{\Phi_1}\Big)\frac{\Phi_1}{C_1 \|w_{5,0}\|_{\infty}}.
\label{14}
\end{equation}
We summarize these estimates with the following lemma.

\begin{lemma}
If $F(\tilde{J}_0)$ satisfies \eqref{13}, then, for all $k\ge 1, F(\tilde{J}_k)\le A$.
In particular, if $F(\tilde{J}_0)$ satisfies \eqref{14}, then
\[
F(\tilde{J}_k)\le A= \frac{\Phi_1}{C_1 q\|w_{5,0}\|_{\infty}}\qquad\,\mbox{for all $k\ge 1$}.
\]
Furthermore, if $F(\tilde{J}_0)$ satisfies \eqref{14}, then
\[
F(J_k)\le F(\tilde{J}_0)\mathrm{exp}\Big(\frac{C_1q \|w_{5,0}\|_{\infty}}{\Phi_1}+1\Big)
\qquad\,\mbox{for all $k\ge 1$}.
\]
\end{lemma}

\par Next, we need to estimate the amount that the solution ``drifts'' from its original
base point due to the source term $G(U)$.
We use $W(U_0(-\infty))={\lim}_{y\to -\infty}W(U_0(y))$ as our base point.
From our scheme,
\begin{equation}
W(U^h(x,y))=W(U_0^h(x,y))+G(U_0^h(x,y))(x-kh).
\end{equation}
We denote $U_k^{\infty}=\lim_{y\to -\infty}U^h((k+1)h-,y)$ and $U_0(-\infty)=\lim_{y\to -\infty}U_0(y)$.
Then $W(U_{0}^{\infty})=W(U_{0}(-\infty))+G(U_{0}(-\infty))h$
and $W(U_{k+1}^{\infty})=W(U_{k}^{\infty})+G(U_{k}^{\infty})h$ for $k\ge 0$.
Since $Z_0(-\infty)=0$,  $G(U_{0}(-\infty))=0$.
We deduce that $W(U_{k+1}^{\infty})=W(U_{k}^{\infty})=W(U_{0}(-\infty))$.
Therefore, for all $(x,y)\in J_k$, we have
\begin{equation}
\begin{split}
\|W(U^h(x,y))-W(U_0(-\infty))\|
&\le \|W(U^h(x,y))-W(U_k^{\infty})||+||W(U_k^{\infty})-W(U_0(-\infty))\|
                        \\& =\|W(U^h(x,y))-W(U_k^{\infty})\|
                        \\& \le TV(W(U^h(x,\cdot)))
                        \\& \le CF(J_k).
\end{split}
\end{equation}

In summary, we have established the following theorem.

\begin{theorem}
If $\mathrm{TV}\big(W(U_0)\big)$ is sufficiently small, then the fractional-step Glimm scheme generates
the approximate solutions $U^h(x,y)$ which exist in the whole domain $\{x\ge 0, y\in \mathbb{R}\}$
and have uniformly bounded total variation in the $y$--direction.
Moreover, there is a null set $N\subset \Pi_{k=0}^{\infty}(-1,1)$
such that, for each $\theta \in \Pi_{k=0}^{\infty}(-1,1)\setminus N$,
there exists a sequence $h_i\to 0$ so that
\begin{equation}
U_{\theta}=\lim _{h_i\to 0}U_{h_i,\theta}
\end{equation}
is an entropy solution to problem \eqref{d3} for system \eqref{d1}--\eqref{d2},
where the limit is taken in $L_{loc}^{1}(\Omega)$.
Moreover, $U_{\theta}$ has uniformly bounded total variation in the $y$--direction.
\end{theorem}

The proof of the convergence part will be given in Section 6.

\section{Initial-Boundary Value Problem}

\par In this section, we are concerned with reacting supersonic flows past Lischitz curved wedges.
The problem can be formulated as the initial-boundary value problem for
system \eqref{d1}--\eqref{d2} in $\Omega$
with initial data on $\Gamma$:
\begin{equation}
(u,v,p,\rho,Z)|_{x=0}=(u_0,v_0,p_0,\rho_0,Z_0)(y)\equiv U_0(y),\qquad \text{$y \in \mathbb{R}$},
\label{dd3}
\end{equation}
and boundary condition
\begin{equation}
(u,v)\cdot \textbf{n}=0 \qquad \text{on $\Gamma$},\label{dd4}\\
\end{equation}
where
\begin{equation*}
\Omega=\{(x,y)\, :\, y<g(x), x>0\},\quad \Gamma=\{(x,y)\, :\, y=g(x), x>0\},
\end{equation*}
and $\textbf{n}(x\pm)=\frac{(-g'(x\pm),1)}{\sqrt{(g'(x\pm))^2+1}}$ is the outer unit normal vector
to $\Gamma$ at the point $x\pm$ (see Fig. 3).

\begin{center}
\setlength{\unitlength}{1mm}
\begin{picture}(120,70)(-2,-10)
\linethickness{1pt}
\put(30,0){\vector(0,1){50}}
\put(0,30){\vector(1,0){120}}
\put(3,3){\vector(4,1){20}}
\put(3,7){\vector(4,1){20}}
\put(3,11){\vector(4,1){20}}
\put(3,-1){\vector(4,1){20}}
\qbezier(30,30)(40,31)(50,33)
\qbezier(50,33)(90,25)(110,32)
\put(24,26){$O$}
\put(116,26){$x$}
\put(27,50){$y$}
\put(6,20){$U_0$}
\put(34,30){\line(1,2){2}}
\put(38,31){\line(1,2){2}}
\put(42,31){\line(1,2){2}}
\put(46,32){\line(1,2){2}}
\put(50,32){\line(1,2){2}}
\put(54,32){\line(1,2){2}}
\put(58,31){\line(1,2){2}}
\put(62,31){\line(1,2){2}}
\put(66,31){\line(1,2){2}}
\put(70,30){\line(1,2){2}}
\put(74,30){\line(1,2){2}}
\put(78,29){\line(1,2){2}}
\put(82,29){\line(1,2){2}}
\put(86,29){\line(1,2){2}}
\put(90,29){\line(1,2){2}}
\put(94,29){\line(1,2){2}}
\put(98,30){\line(1,2){2}}
\put(100,20){\vector(-2,1){16}}
\put(100,18){$y=g(x)$}
\put(52,10){$\Omega$}
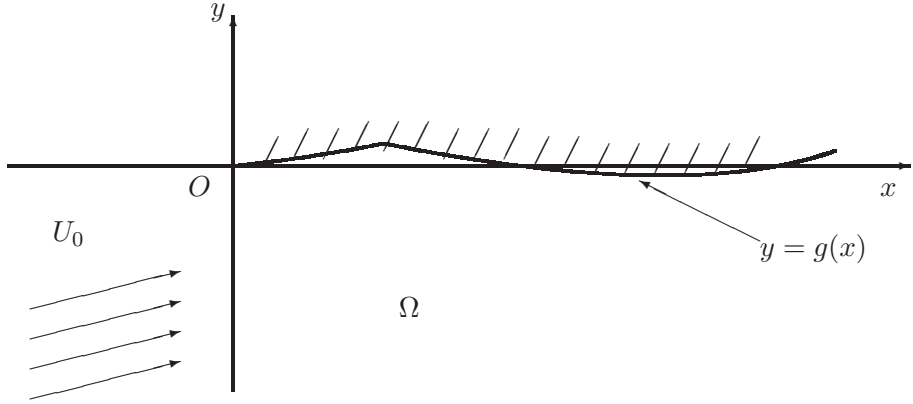
\put(18,-8){Fig. 3. Supersonic flow past a Lipschitz curved wedge}
\end{picture}
\end{center}

The assumption for $U_0(y):= (u_0,v_0,p_0,\rho_0,Z_0)(y)$ is the same as before.
The boundary function $y=g(x)$ is a small perturbation of the straight line $y=\frac{v_0(-\infty)}{u_0(-\infty)}x$
such that $y=g(x)$ is Lipschitz continuous with $g(0)=0, g'(0+)=\arctan(\frac{v_0(-\infty)}{u_0(-\infty)})$,
and $g' \in BV(\mathbb{R}^+;\mathbb{R})$.

Without loss of generality, we may assume that
\begin{equation}
v_0(-\infty)=0, \qquad Z_0(-\infty)=0.
\end{equation}

The formulation of the initial-boundary value problem is derived from the original physical problem
when supersonic flow past a symmetric wedge through the coordinate transformation.
For the non-reacting supersonic flow past a straight symmetric wedge, {\it i.e.} $g'(x)=0$,
a plane shock is generated, which is attached to the wedge vertex (see Fig. 4).
When the supersonic flow is governed by exothermically reacting steady Euler equations,
the attached shock is no longer a plane shock even for the straight wedge,
though it can be handled as an approximate shock wave.

\begin{center}
\setlength{\unitlength}{1mm}
\begin{picture}(100,80)(-2,-10)
\linethickness{1pt}
\put(0,0){\vector(1,0){30}}
\put(-2,10){\vector(1,0){30}}
\put(-4,20){\vector(1,0){30}}
\put(-2,30){\vector(1,0){30}}
\put(0,40){\vector(1,0){30}}
\put(30,20){\line(3,4){20}}
\put(30,20){\line(3,-4){20}}
\put(30,20){\line(5,1){58}}
\put(38,26){\vector(4,1){20}}
\put(30,20){\line(5,-1){58}}
\put(38,13){\vector(4,-1){20}}
\put(41,-4){\it{S}}
\put(43,44){\it{S}}
\put(33,19){\line(5,3){5}}
\put(38,18){\line(5,4){6}}
\put(43,17){\line(4,3){10}}
\put(48,16){\line(5,3){18}}
\put(53,15){\line(5,3){23}}
\put(58,14){\line(5,3){29}}
\put(63,13){\line(5,3){25}}
\put(68,12){\line(5,3){21}}
\put(73,11){\line(5,3){16}}
\put(78,10){\line(5,3){11}}
\put(83,9){\line(5,3){6}}
\put(18,-14){Fig. 4. Non-reacting supersonic flow past a straight wedge}
\end{picture}
\end{center}

\bigskip
\smallskip
\subsection{\textbf{Homogeneous initial-boundary value problem}}

\par We first recall some basic properties on the initial-boundary value problem for the homogeneous system \eqref{3}.

\subsubsection{\textbf{Lateral Riemann problem}}
 The simple case of problem (\ref{d1})--(\ref{d2})
is that $g\equiv 0$.
It has been shown in \cite{courant} that, if $g\equiv 0$, the homogeneous system \eqref{3}
with initial condition:
\begin{equation}
(u,v,p,\rho,Z)|_{x<0}=(u_{-},v_{-},p_{-},\rho_{-},Z_{-})\equiv U_{-}
\end{equation}
yields an entropy solution that consists of the constant states $U_{-}$ and
$U_{+}:=(u_{+},0,p_{+},\rho_{+},Z_{+})$  with $u_{+}>c_{+}>0$ in the subdomain of $\Omega$,
separated by a straight shock-front emanating from the vertex.
That is, the state ahead of the shock-front is $U_{-}$,
whilst the state behind the shock-front is $U_{+}$ (see Figs. 5--6).
When the angle between the flow direction of the front state and the wedge boundary
at a boundary vertex is larger than $\pi$, the entropy solution contains a rarefaction wave
that separates the front state from the back state (see Fig. 6).

\bigskip
\begin{center}
\setlength{\unitlength}{1mm}
\begin{picture}(90,60)(-2,-10)
\linethickness{1pt}
\put(30,0){\vector(0,1){50}}
\put(0,30){\vector(1,0){70}}
\put(3,3){\vector(1,1){20}}
\put(30,30){\line(1,-1){25}}
\put(50,20){\vector(1,0){20}}
\put(30,30){\line(1,2){2}}
\put(34,30){\line(1,2){2}}
\put(38,30){\line(1,2){2}}
\put(42,30){\line(1,2){2}}
\put(46,30){\line(1,2){2}}
\put(50,30){\line(1,2){2}}
\put(54,30){\line(1,2){2}}
\put(58,30){\line(1,2){2}}
\put(62,30){\line(1,2){2}}
\put(22,26){$O$}
\put(66,26){$x$}
\put(23,50){$y$}
\put(6,20){$U_{-}$}
\put(54,14){$U^{+}$}
\put(52,2){Shock}
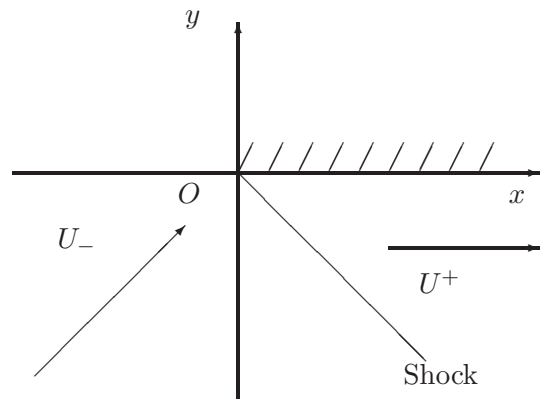
\put(-2,-8){Fig. 5. Unperturbed case when $g\equiv 0$}
\end{picture}
\end{center}

\begin{center}
\setlength{\unitlength}{1mm}
\begin{picture}(180,40)(-2,-4)
\linethickness{1pt}
\put(0,24){\line(1,0){30}}
\multiput(30,24)(3,0){14}{\line(1,0){2}}
\put(2,24){\line(1,2){2}}
\put(6,24){\line(1,2){2}}
\put(10,24){\line(1,2){2}}
\put(14,24){\line(1,2){2}}
\put(18,24){\line(1,2){2}}
\put(22,24){\line(1,2){2}}
\put(26,24){\line(1,2){2}}
\put(30,24){\line(1,2){2}}
\put(34,23){\line(1,1){4}}
\put(38,21){\line(1,1){4}}
\put(42,20){\line(1,1){4}}
\put(46,19){\line(1,1){4}}
\put(50,17){\line(1,1){4}}
\put(54,16){\line(1,1){4}}
\put(58,15){\line(1,1){4}}
\put(62,13){\line(1,1){4}}
\put(70,20){$x$}
\put(68,24){\vector(1,0){3}}
\put(85,24){\line(1,0){30}}
\put(89,24){\line(0,1){3}}
\put(93,24){\line(0,1){3}}
\put(97,24){\line(0,1){3}}
\put(101,24){\line(0,1){3}}
\put(105,24){\line(0,1){3}}
\put(109,24){\line(0,1){3}}
\put(113,24){\line(0,1){3}}
\put(116,24){\line(1,2){2}}
\put(120,25){\line(1,2){2}}
\put(124,26){\line(1,2){2}}
\put(128,27){\line(1,2){2}}
\put(132,28){\line(1,2){2}}
\put(136,29){\line(1,2){2}}
\put(140,30){\line(1,2){2}}
\put(144,31){\line(1,2){2}}
\put(148,32){\line(1,2){2}}
\put(152,33){\line(1,2){2}}
\put(30,24){\line(3,-1){40}}
\put(42,17){\vector(3,-1){18}}
\put(30,24){\line(4,-3){25}}
\put(50,2){Shock}
\put(3,16){\vector(1,0){20}}
\put(115,24){\line(4,1){40}}
\put(126,22){\vector(4,1){20}}
\multiput(115,24)(3,0){13}{\line(1,0){2}}
\put(153,24){\vector(1,0){3}}
\put(115,24){\line(2,-1){18}}
\put(115,24){\line(3,-4){12}}
\put(115,24){\line(1,-2){10}}
\put(90,20){\vector(1,0){16}}
\put(116,0){$\text{Rarefaction wave}$}
\put(154,20){$x$}
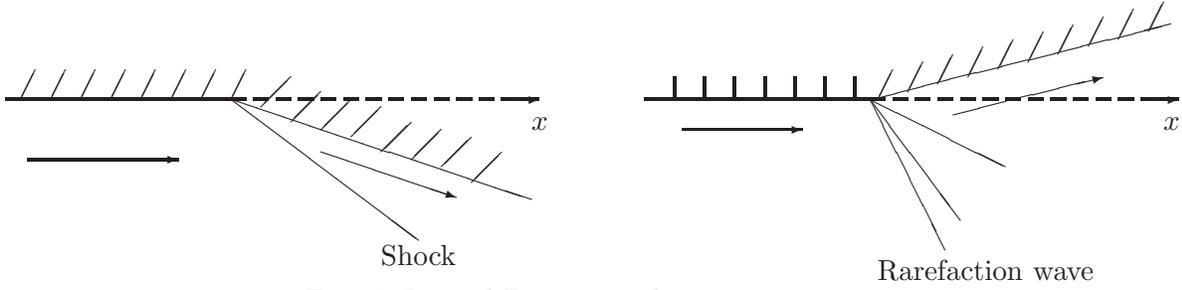
\put(40,-4){Fig. 6. Lateral Riemann solutions}
\end{picture}
\end{center}

\medskip
\subsubsection{\textbf{Riemann problem}}
Consider the Riemann problem for (\ref{3}):
\begin{eqnarray}
U|_{x=x_0}=U_{-}=\left\{
\begin{array}{ll}
U_b, &\quad  y< y_0,\\[2mm]
U_a, &\quad y> y_0,
\end{array}
\right.
\end{eqnarray}
where $U_a$ and $U_b$ are the constant states which are regarded as the above state and below state
with respect to the line $y=y_0$, respectively.
It is well known that this Riemann problem is solvable if the states $U_b$ and $U_a$ are close enough.

\subsubsection{\textbf{Estimates on wave interactions for \eqref{3}}}
The estimates on week wave interactions are the same as in Lemma \ref{4.1b}.

\subsection{\textbf{Estimates of the reflection on the boundary for system \eqref{3}}}
\par Following the notation in \cite{zhang1},
we denote $\{C_k(a_k,b_k)\}_{k=0}^{\infty}$ by the points $\{(a_k,b_k)\}_{k=0}^{\infty}$
in the $(x, y)$--plane with $a_{k+1}>a_k\ge 0$.
Set
\begin{eqnarray}
\omega_{k,k+1}=\arctan\big(\frac{b_{k+1}-b_k}{a_{k+1}-a_k}\big),
\quad \omega_k=\omega_{k,k+1}-\omega_{k-1,k},\quad \omega_{-1,0}=0,\\[2mm]
\Omega_{k}=\{(x,y):x\in [a_k,a_{k+1}), y<b_k+(x-a_k)\tan(\omega_{k,k+1})\},\\[2mm]
\Gamma_{k}=\{(x,y):x\in [a_k,a_{k+1}), y=b_k+(x-a_k)\tan(\omega_{k,k+1})\},
\end{eqnarray}
and the outer unit normal vector to $\Gamma_{k+1}$:
\begin{equation}
\textbf{n}_{k+1}=\frac{(-b_{k+1}+b_k,a_{k+1}-a_k)}{\sqrt{(b_{k+1}-b_k)^2+(a_{k+1}-a_k)^2}}=(-\sin (\omega_{k,k+1}),\cos(\omega_{k,k+1})).
\end{equation}
\par We then consider the initial-boundary value problem:
\begin{eqnarray}
\left\{
\begin{array}{ll}
(\ref{3}) \qquad \text{in $\Omega_{k}$},\\[2mm]
U|_{x=a_k}=\underline{U},\\[2mm]
(u,v)\cdot \textbf{n}_{k}=0 \qquad \text{on $\Gamma_{k}$},
\end{array} \right.
\end{eqnarray}
where $\underline{U}$ is a constant state.

\bigskip
\medskip
\begin{center}
\setlength{\unitlength}{1mm}
\begin{picture}(50,72)(-2,-6)
\linethickness{1pt}
\multiput(0,0)(0,3){21}{\line(0,2){2}}
\multiput(24,2)(0,3){17}{\line(0,1){2}}
\multiput(48,-1)(0,3){20}{\line(0,2){2}}
\put(0,18){\line(1,1){20}}
\put(0,18){\line(5,3){20}}
\put(0,62){\line(5,-2){32}}
\put(23,53){\line(5,-3){18}}
\put(12,57){\vector(1,3){4}}
\put(24,52){\line(4,1){24}}
\put(0,62){\line(1,-1){18}}
\put(36,55){\vector(-1,4){3}}
\put(3,61){\line(5,2){4}}
\put(7,59){\line(5,2){4}}
\put(12,57){\line(5,2){4}}
\put(15,56){\line(5,2){4}}
\put(19,54){\line(5,2){4}}
\put(24,52){\line(1,2){2}}
\put(28,53){\line(1,2){2}}
\put(32,54){\line(1,2){2}}
\put(36,55){\line(1,2){2}}
\put(40,56){\line(1,2){2}}
\put(44,57){\line(1,2){2}}
\put(-1,65){$\Gamma_k$}
\put(49,60){$\Gamma_{k+1}$}
\put(11,70){$n_k$}
\put(35,68){$n_{k+1}$}
\put(9,54){$U_k$}
\put(38,49){$U_{k+1}$}
\put(16,50){$C_{k+1}$}
\put(2,34){$U_m$}
\put(2,10){$U_b$}
\put(36,20){$U_b$}
\put(29,51){$\omega_{k+1}$}
\put(12,4){$\Omega_k$}
\put(36,4){$\Omega_{k+1}$}
\put(36,40){$\delta_1$}
\put(12,43){$\beta_1$}
\put(12,36){$\alpha_5$}
\put(12,23){$\alpha_{2(3,4)}$}
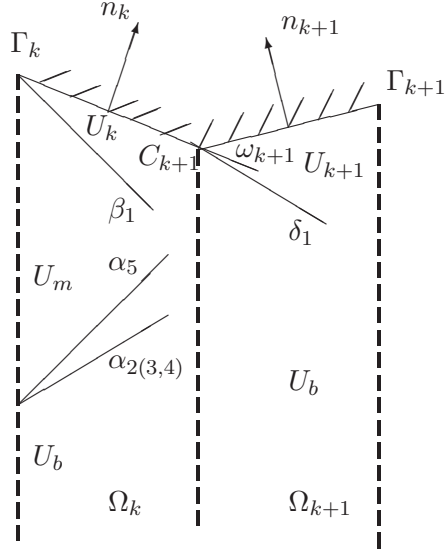
\put(-12,-6){Fig. 8. Weak wave reflections on the boundary.}
\end{picture}
\end{center}

\medskip
\begin{lemma}
Let $\{U_b, U_m\}=(0,\alpha_2,\alpha_3,\alpha_4,\alpha_5)$ and $\{U_m,U_k\}=(\beta_1, 0,0,0,0)$ with
\begin{equation*}
(u_k,v_k) \cdot \textbf{n}_k=0.
\end{equation*}
Then there exists $U_{k+1}$ such that
\begin{equation*}
\{U_b, U_{k+1}\}=(\delta_1,0,0,0,0)\qquad \text{with} \quad (u_{k+1},v_{k+1}) \cdot \textbf{n}_{k+1}=0.
\end{equation*}
Furthermore,
\begin{equation*}
\delta_1=\beta_1+K_{b5}\alpha_5+K_{b4}\alpha_4+K_{b3}\alpha_3+K_{b2}\alpha_2+K_{b0}\omega_k,
\end{equation*}
where $K_{b5},K_{b4},K_{b3},K_{b2}$, and $K_{b0}$ are $C^2$--functions
of $(\alpha_5,\alpha_4,\alpha_3,\alpha_2,\beta_1,\omega_k;U_b)$ satisfying
\begin{equation*}
K_{b5}|_{\omega_k=\alpha_5=\alpha_4=\alpha_3=\alpha_2=\beta_1=0,U_b=U_+}=1,\qquad
K_{bi}|_{\omega_k=\alpha_5=\alpha_4=\alpha_3=\alpha_2=\beta_1=0,U_b=U_+}=0,\,\,\, \text{$i=2,3,4$},
\end{equation*}
and $K_{b0}$ is bounded.
\end{lemma}

\par The proof of this lemma is similar to Proposition 3.2 in \cite{chen2}.

\subsection{\textbf{Construction of approximate solutions}}

\par In this section, we develop a modified Glimm difference scheme to construct a family of
approximate solutions in consistent with the boundary condition (\ref{dd3})--\eqref{dd4}
and establish their necessary estimates for the initial-boundary value
problem for system \eqref{d1}--(\ref{d2})
in the corresponding domains $\Omega_{h}$.

\par We first use the fact that the boundary is a perturbation of the straight wedge:
\begin{equation}
\sup_{x\ge 0}|g'(x)|<\varepsilon \qquad \text{for sufficiently small $\varepsilon>0$.}
\label{612}
\end{equation}
Let $h>0,l>0$ denote the step-length in the $x$-direction and $y$-direction, respectively.
Set $a_k:=kh$ and $b_k:=y_k=g(kh)$ and follow the notations in Section 2.4.
Then
\begin{equation}
m:=\sup_{k>0}\Big\{\frac{|y_k-y_{k-1}|}{h}\Big\}<\varepsilon.
\end{equation}
Define
\begin{equation}
\Omega_{h}=\bigcup_{k\ge 0}\Omega_{h,k},
\end{equation}
where $\Omega_{h,k}=\{(x,y):kh\le x <(k+1)h, \quad y\le g_h(x) \}$
with $g_h(x)=y_k+(x-kh)\tan(\omega_{k,k+1})$ when $kh\le x< (k+1)h$.
We also need the Courant-Friedrichs-Lewy type condition:
\begin{equation}
\max_{1\le j\le 5}\Big(\sup_{U\in O_{\varepsilon}(U_{+})}|\lambda_j(U)|\Big)\le \frac{l-mh}{2h}.
\end{equation}

\par Define
\begin{equation}
a_{k,n}=(2n+1+\theta_k)l+y_k,
\end{equation}
where $\theta_k$ is randomly chosen in $(-1,1)$. Then we choose
\begin{equation}
P_{k,n}=(kh,a_{k,n}),\qquad \text{$k\ge 0, n=0,-1,-2,\cdots $},
\end{equation}
to be the mesh points and define the approximate solutions $W(U^h(x,y))$ in $\Omega_h$
for any $\theta=(\theta_0,\theta_1,\theta_2,\cdots)$ in an inductive way.

\par We denote $T_{k,0}$ the diamond domain whose vertices are $(kh,y_k),(kh,-l+y_k),((k+1)h,-l+y_{k+1})$,
and $((k+1)h,y_{k+1})$.
For $n\le -1$, we denote $T_{k,n}$ the diamond whose vertices are
$(kh,(2n+1)l+y_k),(kh,(2n-1)l+y_k),((k+1)h,(2n-1)l+y_{k+1})$, and $((k+1)h,(2n+1)l+y_{k+1})$.

\par Now we can define the difference scheme in $\Omega_{h}$, that is, define the global approximate
solution $W(U^h(x,y))$ in $\Omega_{h}$. This can be done by carrying out the following steps inductively, similar to
the construction in Section 3.

\par Assume that $W(U^h(x,y))$ is defined for $x<kh$. Then we define $W(U^h(kh+0,y))$ as follows:
\par We define, for $n\le -1$,
\begin{equation}
W(U_0^k):=W(U^h(kh-,a_{k,n})) \qquad \text{for $\,\,\, 2nl+y_k\le y <2(n+1)l+y_k$},
\end{equation}
and
\begin{equation}
W(U^h(kh+0,y)):=W(U_0^k).
\end{equation}

\par First, we define $W(U_0^h(x,y))$ in $T_{k,0}$ by solving the following lateral Riemann problem:
\begin{equation}
\left \{
\begin{array}{ll}
 W(U_k)_{x}+H(U_k)_y=0 \qquad \text{in $T_{k,0}$},\\[2mm]
 W(U_k)|_{x=kh}=W(U_0^k),\\[2mm]
 (u_k,v_k)\cdot \textbf{n}_k=0 \qquad \text{on $\Gamma_k$}.
 \end{array} \right.
\end{equation}
We can obtain the above lateral Riemann solution $W(U_k)$ in $T_{k,0}$ and define
\begin{equation}
W(U_0^h)=W(U_k) \qquad \text{in $T_{k,0}$}.
\end{equation}

\par Second, we solve the following Riemann problem in each diamond $T_{k,n}$ for $n\le -1$:
\begin{equation}
\left \{
 \begin{array}{ll}
   W(U_k)_{x}+H(U_k)_y=0 \qquad \text{in $T_{k,n}$},\\[2mm]
   W(U_k)|_{x=kh}=W(U_0^k),
 \end{array} \right.
\end{equation}
to obtain the Riemann solution $W(U_k)$ in $T_{k,n}$ and define
\begin{equation}
W(U_0^h)=W(U_k) \qquad \text{in $T_{k,n}, n\le -1$}.
\end{equation}

\par Finally, we use the Glimm fractional-step operator to obtain the desired approximate solutions:
\begin{equation}
W(U^h(x,y))=W(U_0^h(x,y))+G(U_0^h(x,y))(x-kh)\qquad \text{for $kh\le x<(k+1)h$}.
\end{equation}

\par In this way, we have constructed the approximate solution $W(U^h(x,y))$ globally,
provided that we can obtain a uniform bound of the approximate solutions.
To achieve this, we establish the total variation of $W(U^h(x,y))$ on a class of space-like curves.

\par As before, for the mesh curves $J$ in $x>0$, we give the following definition:
\begin{definition}
\begin{align*}
 & L_0(J)=\sum\{\omega(C_k):C_k \in \Omega_J\},\\
 & L_j(J)=\sum\{|\alpha_j|:\text{$\alpha_j$ cross $J$}\},\quad j=1,2,3,4,5,\\
 & L(J)=K^*L_0(J)+L_1(J)+K^*\big(L_2(J)+L_3(J)+L_4(J)+L_5(J)\big),\\
 & Q(J)=\sum\{|\alpha_i||\beta_j|:\text{both $\alpha_i$ and $\beta_j$ crossing $J$ and approaching}\},\\
 & F(J)=L(J)+K Q(J),
\end{align*}
where $K>0$ is determined later, $\Omega_J$ is the set of the corner points $C_k$ with $k\ge 0$:

\begin{equation}
\Omega_J=\{C_k \, :\, \,\, C_k \in J\cap \partial \Omega_{h},\mbox{$C_k=(kh,g(kh)),k\ge 0$}\},
\end{equation}
and $K^*$ is a positive constant that satisfies $K^*>\max_{2\le i\le 5} K_{bi}+1$.
\end{definition}

\par Next, we estimate the functional $F$.
To do this, let $I$ and $J$ be two $k$-mesh curves for some $k>0$ such that $J$ is an immediate
successor to $I$, and let $\Lambda$ be the diamond between $I$ and $J$.
Due to the location of $\Lambda$, two cases are to be considered:

\begin{enumerate}
\item[\rm (i)] Case $\Lambda \subset \Omega_h$: If $\alpha$ and $\beta$ are the waves entering $\Lambda$,
we define
\begin{equation}
Q(\Lambda)=\sum|\alpha_i||\beta_j|,
\end{equation}
where the sum is taken over all the pairs for which the $i$-wave from $\alpha$
and $j$-wave from $\beta$ are approaching;

\medskip
\item[\rm (ii)] Case $\Lambda \cap \partial \Omega_h \neq \emptyset$:
Let $\Omega_J=\Omega_I \setminus \{C_k\}$ with $C_k=(kh,y_k)$ for some $k\ge 0$,
let $I=I_0\cup I'$ and $J=I_0\cup J'$ such that $\partial \Lambda=I'\cup J'$, and
let $\beta_1$ and $\alpha_i$ be the $1$-wave and $i$-wave respectively crossing
$I'$ with $\alpha_i$ lying below $\beta_1$ on $I$, where $i=2,3,4,5$.
In addition, by the construction of approximate solutions,
let $\delta_1$ be the weak $1$-wave crossing $J'$ (see Fig. 9 below).
\end{enumerate}

\par Define
\begin{eqnarray}
E_{h,\theta}(\Lambda)=\left\{
\begin{array}{ll}
\omega_k+\sum_{i=2}^{5}|\alpha_i| &\quad \text{if $\Lambda \cap \partial \Omega_h \neq \emptyset$},\\[2mm]
Q(\Lambda) & \quad \text{if $\Lambda \subset \Omega_h$}.
\end{array}
\right.
\end{eqnarray}

\subsection{\textbf{Estimates of the non-reacting step involving the boundary}}

\par By choosing a suitable constant $K$,
we now prove that the Glimm-type functional $F$ is non-increasing in the non-reacting step.

\begin{theorem}
Suppose that the wedge function $g(x)$ satisfies {\rm (\ref{612})}, and $I$ and $J$ are two mesh curves
such that $J$ is an immediate successor of $I$. Then there exist constants $\varepsilon>0$ and $K>0$
such that, if $F(I)\le \varepsilon$, then
\begin{equation}
F(J)\le F(I)-\frac{1}{4}E_{h,\theta}(\Lambda).
\end{equation}
\end{theorem}

\begin{pf}
{We divide our proof into two cases depending on the location of the diamond.}

\medskip
\par {\textbf{Case 1} (interior weak-weak interaction):
$\Lambda$ lies in the interior of $\Omega_h$. Denote $Q(\Lambda)=\Delta(\alpha,\beta)$ as
defined in Lemma {\rm 2.1}. Then, for some constant $M>0$,}
\begin{equation}
L(J)-L(I)\le (1+4K^*)MQ(\Lambda).
\end{equation}
Since $L(I_0)<\varepsilon$ from $F(I)<\varepsilon$, we have
\begin{equation}
\begin{split}
Q(J)-Q(I)&=\big(Q(I_0)+\sum_{i=1}^{5}Q(\gamma_i,I_0)\big)-\big(Q(I_0)+Q(\Lambda)+\sum_{i=1}^{5}Q(\alpha_i,I_0)+
\sum_{i=1}^{5}Q(\beta_i,I_0)\big)
\\&\le Q(MQ(\Lambda),I_0)-Q(\Lambda)\\[1.5mm]
&\le \big(ML(I_0)-1\big)Q(\Lambda)\\
&\le -\frac{1}{2}Q(\Lambda).
\end{split}
\end{equation}
Hence, by choosing a suitably large $K$, we obtain
\begin{equation}
F(J)-F(I)\le \big((1+4K^*)M-\frac{K}{2}\big)Q(\Lambda)\le -\frac{1}{4}Q(\Lambda).
\end{equation}

\par \textbf{Case 2} (near the boundary):
$\Lambda$ touches the approximate boundary $\partial \Omega_h$.
Then $\Omega_J=\Omega_I \setminus \{C_k\}$ for certain $k$.

\par Let $\delta_1$ be the weak $1$-wave going out of $\Lambda$ through $J'$,
and let $\beta_1,\alpha_2,\alpha_3,\alpha_4$, and $\alpha_5$ be the weak waves
entering $\Lambda$ through $I'$, as shown in Fig. 9. Then

\begin{center}
\setlength{\unitlength}{1mm}
\begin{picture}(50,76)(-2,0)
\linethickness{1pt}
\multiput(0,0)(0,3){21}{\line(0,2){2}}
\multiput(24,2)(0,3){22}{\line(0,2){2}}
\multiput(48,-1)(0,3){20}{\line(0,2){2}}
\put(0,18){\line(1,1){20}}
\put(0,18){\line(5,3){21}}
\put(24,26){\line(4,3){24}}
\put(24,26){\line(3,-4){12}}
\put(0,62){\line(5,-2){32}}
\put(23,53){\line(5,-4){20}}
\put(24,52){\line(4,1){24}}
\put(0,62){\line(1,-1){18}}
\put(24,66){\line(-1,-1){24}}
\put(24,65){\line(6,-5){24}}
\put(0,42){\line(3,-2){24}}
\put(3,61){\line(5,2){4}}
\put(7,59){\line(5,2){4}}
\put(12,57){\line(5,2){4}}
\put(15,56){\line(5,2){4}}
\put(19,54){\line(5,2){4}}
\put(24,52){\line(1,2){2}}
\put(28,53){\line(1,2){2}}
\put(32,54){\line(1,2){2}}
\put(36,55){\line(1,2){2}}
\put(40,56){\line(1,2){2}}
\put(44,57){\line(1,2){2}}
\put(19,49){$C_k$}
\put(37,10){$I_0$}
\put(2,35){$I'$}
\put(36,32){$J'$}
\put(29,51){$\omega_k$}
\put(36,44){$\delta_1$}
\put(12,43){$\beta_1$}
\put(12,36){$\alpha_5$}
\put(12,23){$\alpha_{2(3,4)}$}
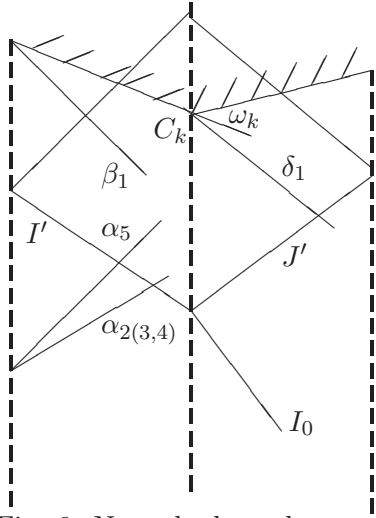
\put(-2,-4){Fig. 9. Near the boundary.}
\end{picture}
\end{center}

\begin{equation*}
L_0(J)-L_0(I)=-|\omega_k|,
\end{equation*}
\begin{equation*}
L_i(J)-L_i(I)=\sum_{\text{$\gamma_i$ cross $I_0$}}|\gamma_i|
 -\big(|\alpha_i|+\sum_{\text{$\gamma_i$ cross $I_0$}}|\gamma_i|\big)
 =-|\alpha_i|, \qquad \text{$i=2,3,4,5$},
\end{equation*}
\begin{equation*}
\begin{split}
L_1(J)-L_1(I)&=\Big(|\delta_1|+\sum_{\text{$\gamma_1$ cross $I_0$}}|\gamma_1|\Big)
  -\Big(|\beta_1|+\sum_{\text{$\gamma_1$ cross $I_0$}}|\gamma_1|\Big)
\\& =|\delta_1|-|\beta_1|\\
&\le \sum_{i=2}^{5}|K_{bi}||\alpha_i|+|K_{b0}||\omega_k|,
\end{split}
\end{equation*}
where the last step is from Lemma 5.1. Thus,
\begin{equation}
\begin{split}
L(J)-L(I)&\le (|K_{b0}|-K^*)|\omega_k|+\sum_{i=2}^{5}(|K_{bi}|-K^*)|\alpha_i|
\\& \le-\big(|\omega_k+\sum_{i=2}^{5}|\alpha_i|\big),
\end{split}
\end{equation}
since $K^*>\max K_{bi}+1$ for $ i=1,2,3,4,5$.  Moreover, we have
\begin{equation}
\begin{split}
Q(J)-Q(I)&=\big(Q(I_0)+Q(\delta_1,I_0)\big)
-\Big(Q(I_0)+Q(\beta_1,I_0)+\sum_{i=2}^{5}Q(\alpha_i,I_0)+|\beta_1|\sum_{i=2}^{5}|\alpha_i|\Big)
\\& \le \Big(\sum_{i=2}^{5}|K_{bi}||\alpha_i|+|K_{b0}||\omega_k|\Big)L(I_0).
\end{split}
\end{equation}
Then we obtain
\begin{equation}
\begin{split}
F(J)-F(I)&=\big(L(J)-L(I)\big)+K\big(Q(J)-Q(I)\big)
\\& \le -\Big(|\omega_k|+\sum_{i=2}^{5}|\alpha_i|\Big)
+K\Big(\sum_{i=2}^{5}|K_{bi}||\alpha_i|+|K_{b0}||\omega_k|\Big)L(I_0)
\\& \le -\frac{1}{4}\Big(|\omega_k|+\sum_{i=2}^{5}|\alpha_i|\Big),
\end{split}
\end{equation}
since we can choose $\varepsilon$ sufficiently small.
The proof is completed.
\end{pf}

\subsection{\textbf{Estimates of the reacting step involving the boundary}}

\par We first consider the change of the wave strength before and after reaction
near the boundary. We denote by $(\tilde{U_b}, \tilde{U}_*)$ and $\tilde{\beta_1}$
the two states and wave strength before reaction, respectively,
while by $(U_b, U_*)$ and $\beta_1$ after reaction, respectively (see Figure 10).
According to the boundary condition, we have
\begin{center}
\setlength{\unitlength}{1mm}
\begin{picture}(100,78)(-2,-6)
\linethickness{1pt}
\multiput(0,0)(0,3){13}{\line(0,2){2}}
\multiput(25,0)(0,3){15}{\line(0,2){2}}
\multiput(73,-1)(0,3){13}{\line(0,2){2}}
\multiput(98,-1)(0,3){15}{\line(0,2){2}}
\put(0,38){\line(4,1){25}}
\put(73,37){\line(4,1){25}}
\put(73,37){\line(4,-1){20}}
\put(0,38){\line(4,-1){20}}
\put(21,33){$\tilde{\beta_1}$}
\put(14,37){$\tilde{U}_*$}
\put(94,33){$\beta_1$}
\put(87,37){$U_*$}
\put(10,26){$\tilde{U_b}$}
\put(83,26){$U_b$}
\put(15,42){\vector(-1,3){4}}
\put(88,41){\vector(-1,3){4}}
\put(15,50){$n_{k}$}
\put(88,50){$n_{k}$}
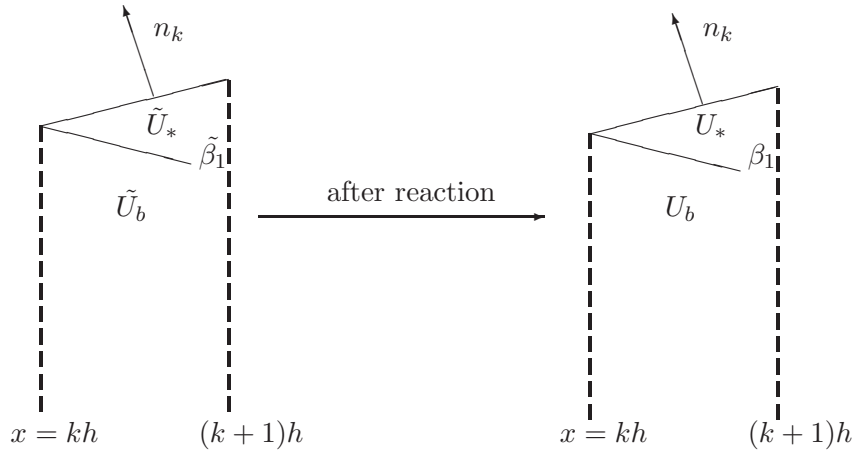
\put(4,-11){Fig. 10. Change of wave strength near the boundary}
\put(38,28){after reaction}
\put(29,26){\vector(1,0){38}}
\put(-4,-4){$x=kh$}
\put(69,-4){$x=kh$}
\put(21,-4){$(k+1)h$}
\put(94,-4){$(k+1)h$}
\end{picture}
\end{center}

\bigskip
\noindent
where
\begin{equation}
W(U_b(x,y))=W(\tilde{U}_b(x,y))+G(\tilde{U}_b(x,y))(x-kh),\qquad \text{$kh\le x<(k+1)h$},
\end{equation}
and
\begin{equation}
W(U_*(x,y))=W(\tilde{U}_*(x,y))+G(\tilde{U}_*(x,y))(x-kh),\qquad \text{$kh\le x<(k+1)h$}.
\end{equation}
From Lemma \ref{q1}, $U_b-\tilde{U}_b=\|Z_0\|_{\infty}e^{-\Phi_1 kh}O(h)$
and $U_*-\tilde{U}_*=\|Z_0\|_{\infty}O(h)e^{-\Phi_1 kh}$.
Therefore, we obtain
\begin{equation}
\beta_1-\tilde{\beta}_1=\|Z_0\|_{\infty}O(h)e^{-\Phi_1 kh}.\label{key}
\end{equation}

As to the inner part, if we perform the same procedure as in the case
of the Cauchy problem, we can obtain a similar estimate:
\begin{equation}
L(J_k)- L(\tilde{J_k})\le Ch\|w_{5,0}\|_{\infty}e^{-\Phi_1 kh}L(\tilde{J_k}).
\end{equation}
Combining these two parts together, we have the following global estimate:
\begin{equation}
L(J_k)- L(\tilde{J_k})\le Ch\|w_{5,0}\|_{\infty}e^{-\Phi_1 kh}\big(L(\tilde{J_k})+1\big).
\end{equation}
Therefore, we can do the same procedure as before to establish

\begin{theorem}
If $\mathrm{TV}\big(W(U_0)\big)+\mathrm{TV}(g')$ is sufficiently small,
then the fractional-step Glimm scheme generates
the approximate solutions $U^h(x,y)$ which exist in the whole domain $\Omega$
and have uniformly bounded total variation in the $y$--direction.
Moreover, there is a null set $N\subset \Pi_{k=0}^{\infty}(-1,1)$
such that, for each $\theta \in \Pi_{k=0}^{\infty}(-1,1)\setminus N$, there exist a sequence $h_i\to 0$
so that
\begin{equation}
U_{\theta}=\lim _{h_i\to 0}U_{h_i,\theta}
\end{equation}
is a weak solution to problem  \eqref{dd3}--\eqref{dd4} for system \eqref{d1}--\eqref{d2},
where the limit is taken in $L_{loc}^{1}(\Omega)$.
Moreover, $U_{\theta}$ has uniformly bounded total variation in the $y$--direction.
\end{theorem}

The proof of the convergence part of Theorem 5.2 is in Section 6.

\section{\textbf{Convergence to Entropy Solutions}}
In this section we show that the limit function of the approximate solutions is an entropy solution
to the Cauchy problem \eqref{d3}--\eqref{d3-a} and
the initial-boundary value problem \eqref{dd3}--\eqref{dd4} for system \eqref{d1}--\eqref{d2}.

\par Let $d\theta_k$ denote the uniform probability measure on $(-1,1)$,
and let $d\theta$ denote the induced product probability measure for the random
sample $\{\theta_k\}_{k=1}^{\infty}$ in the Cartesian product space $\mathscr{A}=\prod_{k=1}^{\infty}(-1,1)$.

\begin{theorem}
Suppose that

\begin{enumerate}
\item[\rm (i)]
The sequence $U^h(x,y)$ is constructed by using the Glimm fractional-step scheme
with the random sample $\{\theta_k\}_{k=0}^{\infty}$ chosen from $\mathscr{A}$.

\item[\rm (ii)]
There exist a null set $\mathscr{N}\subset \mathscr{A}$ such that,
for $\{\theta_k\}\subset \mathscr{A}-\mathscr{N}$,
the sequence $U^h(x,y)$ is uniformly bounded in $L^{\infty}$ and
converges pointwise a.e. to the function $U(x,y)$.
\end{enumerate}
Then the function $U(x,y)$ is an entropy solution of the corresponding problem \eqref{d3}--\eqref{d3-a},
or problem \eqref{dd3}--\eqref{dd4},
for system \eqref{d1}--\eqref{d2}.
That is, for any convex entropy pair $(\eta,q)$ with respect to $W(U)$, the following inequality
\begin{equation}
\eta(W(U))_x+q(W(U))_y\le \nabla_W \eta(W(U))G(U)
\label{u}
\end{equation}
holds in the sense of distributions in $\mathbb{R}^2$ for problem \eqref{d3}--\eqref{d3-a}
and in $\Omega$ including the boundary for problem \eqref{dd3}--\eqref{dd4},
which means that
\begin{eqnarray}
\iint\limits_{\Omega}\big((\eta(W(U))\phi_x+q(W(U))\phi_y+\nabla_W \eta(W(U))G(U)\phi\big)dxdy
\\+\int_{-\infty}^{\infty}\eta(W(U_0(y)))\phi(0,y)dy\ge 0,
\end{eqnarray}
where $\phi(x,y)\ge 0$:
for the Cauchy problem \eqref{d3} with $\Omega=\mathbb{R}^2$ and $\phi \in C_0^{\infty}(\mathbb{R}^2)$;
and  for the initial-boundary value problem \eqref{dd3}--\eqref{dd4},
either $\phi \in C_0^{\infty}(\Omega)$, or
 $\phi \in C_0^{\infty}(\mathbb{R}^2)$ and $(\eta, q)=\alpha(W(U))(u,v)$ for
 any smooth function $\alpha(W)$ of $W$.
\end{theorem}

\begin{pf}  We focus our proof on the initial-boundary value problem \eqref{d3}--\eqref{d3-a},
since the proof for the Cauchy problem \eqref{dd3}--\eqref{dd4} is simpler.

We define
\begin{equation}
\begin{split}
L(\theta,h,\phi) =&\iint\limits_{\Omega_h}\big(\eta(W(U^h))\phi_x+q(W(U^h))\phi_y+\nabla_W \eta(W(U^h))G(U^h)\phi\big)dxdy\\
& +\int_{-\infty}^{0}\eta(W(U_0(y)))\phi(0,y)dy.
\end{split}\label{54}
\end{equation}
We only need to prove that $\lim\limits_{h\to 0} L(\theta,h,\phi)\ge 0$ for $\{\theta_k\}\subset \mathscr{A}-\mathscr{N}$.

\par Since $U_0^h(x,y)$ is an entropy solution of conservation laws $W(U)_x+H(U)_y=0$ in the domain $\Omega_{h,k}$, then
\begin{equation}
\begin{split}
&\iint\limits_{\Omega_{h,k}}\big(\eta(W(U_0^h))\phi_x+q(W(U_0^h))\phi_y\big)dxdy
+\int_{-\infty}^{y_k}\eta(W(U_0^h(kh+0,y)))\phi(kh,y)\, dy\\
&-\int_{-\infty}^{y_{k+1}}\eta(W(U_0^h((k+1)h-,y)))\phi((k+1)h-,y)\,dy\ge 0,
\end{split} \label{55}
\end{equation}
that is,
\begin{equation}
\begin{split}
&\iint\limits_{\Omega_{h,k}}\big(\eta(W(U_0^h))\phi_x+q(W(U_0^h))\phi_y\big)dxdy
+\int_{-\infty}^{0}\eta(W(U_0^h(kh+0,y+y_k)))\phi(kh,y+y_k)\,dy\\
&-\int_{-\infty}^{0}\eta(W(U_0^h((k+1)h-,y+y_{k+1})))\phi((k+1)h-,y+y_{k+1})\,dy\ge 0.
\end{split}\label{55}
\end{equation}
Here we have used the fact that $(u_0^h,v_0^h)\cdot n_k=0$ on the boundary, and the assumptions
for $(\eta,q)$ and $\phi$.
Since $W(U^h(x,y))=W(U_0^h(x,y))+G(U_0^h(x,y))(x-kh)$, then
\begin{equation}
\begin{split}
&\eta(W(U^h(x,y)))-\eta(W(U_0^h(x,y)))\\
&=\nabla_W \eta(W(U_0^h(x,y)))G(U_0^h(x,y))(x-kh)+\varepsilon(x-kh;x,y)(x-kh)
\end{split}
\end{equation}
for some function $\varepsilon(s;x,y)$, which converges uniformly to $0$ as $s\to 0$.
Multiplying the above equation by $\phi_x$ on both sides and integrating on $\Omega_{h,k}$,
we have
\begin{equation}
\begin{split}
&\iint\limits_{\Omega_{h,k}}\big(\eta(W(U^h))-\eta(W(U_0^h))\big)\phi_xdxdy\\
&=\iint\limits_{\Omega_{h,k}}\varepsilon(x-kh;x,y)(x-kh)\phi_xdxdy
 +\iint\limits_{\Omega_{h,k}}\nabla_W \eta(W(U_0^h(x,y)))G(U_0^h(x,y))(x-kh)\phi_xdxdy\\
&=\iint\limits_{\Omega_{h,k}}\varepsilon(x-kh;x,y)(x-kh)\phi_xdxdy
 -\iint\limits_{\Omega_{h,k}}\nabla_W \eta(W(U_0^h(x,y)))G(U_0^h(x,y))\phi_xdxdy\\
&\quad -\iint\limits_{\Omega_{h,k}}\frac{\partial}{\partial x}\big(\nabla_W \eta(W(U_0^h(x,y)))G(U_0^h(x,y))\big)(x-kh)\phi\, dxdy\\
&\quad +\int\limits_{\Gamma_k}\nabla_W \eta(W(U_0^h(x,y)))G(U_0^h(x,y))(x-kh)\phi n_k^1ds\\
&\quad +h\int_{-\infty}^{0}\nabla_W\eta(W(U_0^h((k+1)h-,y+y_{k+1})))G(U_0^h(k+1)h-,y+y_{k+1}))\\
&\qquad\qquad\quad\times \phi(k+1)h-,y+y_{k+1})\,dy.
\end{split}
\end{equation}
Therefore, we use equation (\ref{55}) to obtain
\begin{equation}
\begin{split}
&\iint\limits_{\Omega_{h,k}}\eta(W(U^h))\phi_x\\
&\ge -\iint\limits_{\Omega_{h,k}}\big(q(W(U_0^h))\phi_y+\nabla_W \eta(W(U_0^h(x,y)))G(U_0^h(x,y))\phi\big) dxdy\\
&\quad +\iint\limits_{\Omega_{h,k}}\varepsilon(x-kh;x,y)(x-kh)\phi_x\, dxdy\\
&\quad +\int_{-\infty}^{0}\big(\eta(W(U_0^h((k+1)h-,y+y_{k+1})))\phi((k+1)h,y+y_{k+1})\\
&\qquad\qquad\,\,\,\,\,  -\eta(W(U_0^h(kh+0,y+y_k)))\phi(kh,y+y_k)\big)dy\\
&\quad -\iint\limits_{\Omega_{h,k}}\nabla_W \eta(W(U_0^h(x,y)))G(U_0^h(x,y))\phi\, dxdy\\
&\quad -\iint\limits_{\Omega_{h,k}}\frac{\partial}{\partial x}\big(\nabla_W \eta(W(U_0^h(x,y)))G(U_0^h(x,y))\big)(x-kh)\phi\, dxdy\\
&\quad +\int\limits_{\Gamma_k}\nabla_W \eta(W(U_0^h(x,y)))G(U_0^h(x,y))(x-kh)\phi n_k^1\, ds\\
&\quad +h\int_{-\infty}^{0}\nabla_W\eta(W(U_0^h((k+1)h-,y+y_{k+1})))G(U_0^h(k+1)h-,y+y_{k+1}))\\
&\qquad\qquad\quad\times \phi(k+1)h-,y+y_{k+1})\, dy.
\end{split}
\end{equation}

Summing over $k$, we have
\begin{equation}
L(\theta,h,\phi)\ge \mathscr{A}(\theta,h,\phi)+\sum_{k=0}^{\infty}\mathscr{B}_k(\theta,h,\phi)
+\sum_{k=0}^{\infty}\mathscr{C}_k(\theta,h,\phi)
+\sum_{k=0}^{\infty}\mathscr{D}_k(\theta,h,\phi),
\end{equation}
where
\begin{equation*}
\begin{split}
\mathscr{A}(\theta,h,\phi)=&\sum_{k=0}^{\infty}\mathscr{A}_k(\theta,h,\phi),\\
\mathscr{A}_0(\theta,h,\phi)=&\int_{-\infty}^{0}\big(\eta(W(U_0(y)))-\eta(W(U_0^h(0,y)))\big)\phi(0,y)dy,\\
\mathscr{A}_k(\theta,h,\phi)=&\int_{-\infty}^{0}
\big(\eta(W(U_0^h(kh-,y+y_k)))-\eta(W(U_0^h(kh+0,y+y_k)))\big)\phi(kh,y+y_k)dy\\
&\, + h\int_{-\infty}^{0}\nabla_W\eta(W(U_0^h((k+1)h-,y+y_{k+1})))\times\\
&\, \qquad\qquad \times G(U_0^h(k+1)h-,y+y_{k+1}))\phi(k+1)h-,y+y_{k+1})dy,\qquad k\ge 1,\\
%
\mathscr{B}_k(\theta,h,\phi)=&\iint\limits_{\Omega_{h,k}}\big(q(W(U^h))-q(W(U_0^h))\phi_y\big) dxdy\\
&\, +\iint\limits_{\Omega_{h,k}}\big(\nabla_W \eta(W(U^h))G(U^h)-\nabla_W \eta(W(U_0^h))G(U_0^h)\big)\phi dxdy \\
&\, +\iint\limits_{\Omega_{h,k}}\varepsilon(x-kh;x,y)(x-kh)\phi_xdxdy,\\
\mathscr{C}_k(\theta,h,\phi)=&\iint\limits_{\Omega_{h,k}}\frac{\partial}{\partial x}
\big(\nabla_W \eta(W(U_0^h(x,y)))G(U_0^h(x,y))\big)(x-kh)\phi dxdy,\\
\mathscr{D}_k(\theta,h,\phi)=&\int\limits_{\Gamma_k}\nabla_W \eta(W(U_0^h(x,y)))G(U_0^h(x,y))(x-kh)\phi n_k^1ds.
\end{split}
\end{equation*}
The proof for each component converging to zero as $h$ tends to zero is similar to \cite{chen1},
and we omit here.
This completes the proof.
\end{pf}

\section{\textbf{Asymptotic Behavior involving the Boundary}}

\par Let $\theta \in \Pi_{k=0}^{\infty}(-1,1)\setminus \mathcal{N}$ be equidistributed.
To determine the asymptotic behavior of the solution $U(x,y)$, we need further estimates on $U_{h,\theta}$.

\begin{lemma}
There exists a constant $M_1>0$, independent of $U_{h,\theta}, \theta$ and $h$, such that
\begin{equation}
\sum_{\Lambda}E_{h,\theta}(\Lambda)\le M_1,
\end{equation}
where the summation is over all the diamonds.
\end{lemma}

\begin{pf}
\par First, from the conclusion of the non-reacting step, {\it i.e.} Theorem 5.1, we know
\begin{equation}
F(J)-F(I)\le -\frac{1}{4}E_{h,\theta}(\Lambda),
\end{equation}
where $J$ is an immediate successor of $I$.
Then we conclude
\begin{equation}
F(\tilde{J}_k)-F(J_{k-1})\le -\frac{1}{4}\sum_{k-1}^{k+1}E_{h,\theta}(\Lambda),
\end{equation}
where the summation is over all the diamonds between $x=(k-1)h$ and $x=(k+1)h$.

\par Then we know from the reacting step that
\begin{equation}
F(J_k)-F(\tilde{J_k})\le Ch\|w_{5,0}\|_{\infty}e^{-\Phi_1 kh}\big(F(\tilde{J_k})+1\big)^2.
\end{equation}
Combine these two steps together and sum for $k$ from $1$ to $\infty$ to obtain
\begin{equation*}
\begin{split}
\sum_{k=1}^{\infty}\sum_{k-1}^{k+1}E_{h,\theta}(\Lambda)
& \le CF(J_0)+\sum_{k=1}^{\infty}Ch||w_{5,0}||_{\infty}e^{-\Phi_1 kh}(F(\tilde{J_k})+1)^2\\
& \le C\big(F(J_0)+\|w_{5,0}\|_{\infty}\big)<\infty.
\end{split}
\end{equation*}
The proof is completed.
\end{pf}

\par Moreover, let $\Gamma_g=\cup_{k=0}^{\infty}\bar{\Lambda}_{k,0}$,
where $\Lambda_{k,0}$ is the diamond centered at $C_k$,
and let $L_{h,\theta}(\Gamma_g)$ be the summation of the strength of waves leaving $\Gamma_g$.
Then we have
\begin{lemma}
There exists a constant $M_2$ independent of $U_{h,\theta}, h$, and $\theta$ such that
\begin{equation}
L_{h,\theta}(\Gamma_g)\le M_2 \sum_{\Lambda}E_{h,\theta}(\Lambda).
\end{equation}
\end{lemma}

This can be obtained by employing Lemmas 5.1--5.2 and (\ref{key}) and by taking
the summation of them over $\Gamma_g$.

\par For $i=2, 3, 4, 5$, let $L_i(a-)$  be the amount of all $i$-waves in $U_{\theta}$
crossing the line $x=a$ for any $a>0$.
Also, let $\tilde{L}_{i}^{h,\theta}(a)$ and $L_{i}^{h,\theta}(a)$ denote
the amount of $i$-waves before reaction and after reaction, respectively, in $U_{h,\theta}$
crossing the line $x=a$ for any $a>0$.

\begin{lemma}
$L_i(x-)\to 0$ as $x\to \infty$, for $i=2, 3, 4 ,5$.
\end{lemma}

\begin{pf}
In fact, for $kh\le x< (k+1)h$,
\begin{equation}
 \tilde{L}_{i}^{h,\theta}(x)-L_{i}^{h,\theta}(x)\le L(J_k)-L(\tilde J_k)
 \le  Ch\|w_{5,0}\|_{\infty}e^{-\Phi_1 kh}\big(L(\tilde{J_k})+1\big).
\end{equation}
Then, by Lemmas 7.1--7.2, we can perform the same procedure
as in \cite{zhang2} and conclude this result.
\end{pf}

\par Next, we study the asymptotic behavior of the trace of $U$ on
the boundary.
To this end, from Lemmas 7.1--7.2, we can first deduce

\begin{lemma}
Let
\begin{equation}
B_{h,\theta}(x)=U_{h,\theta}(x, g_h(x)).
\end{equation}
Then there exists a constant $M>0$ depending only on the system such that
\begin{equation}
\mathrm{TV}\{B_{h,\theta}; [0,\infty)\}\le M.
\end{equation}
\end{lemma}

\par Then, by Lemma 7.4,  we can choose a subsequence $\{h_{i_l}\}$ of $\{h_i\}$ so that
\begin{equation}
B_{h_{i_l},\theta} \to B_{\theta}\label{ee}
\end{equation}
in $L_{loc}^1([0,\infty))$ as $h_{i_l} \to 0$ for
some $B_{\theta}\in L^{\infty}$.
From the construction of approximate solutions, we have

\begin{lemma}
Let $B_{\theta}$ be given by \eqref{ee}. Then
\begin{equation*}
B_{\theta}\in BV([0,\infty))
\end{equation*}
and
\begin{equation*}
B_{\theta}(x-)\cdot (-g'(x-),1,0,0,0)=0.
\end{equation*}
\end{lemma}

\begin{pf} Since
\begin{equation}
\begin{split}
&B_{h_{i_l},\theta}(x-)\cdot (-g'_{h_{i_l}}(x-),1,0,0,0)\\
&=(B_{h_{i_l},\theta}(x-)-\tilde{B}_{h_{i_l},\theta}(x-))\cdot (-g'_{h_{i_l}}(x-),1,0,0,0)
+\tilde{B}_{h_{i_l},\theta}(x-)\cdot (-g'_{h_{i_l}}(x-),1,0,0,0),
\end{split}
\end{equation}
the first term on the right-hand side tend to $0$ as $h_{i_l} \to 0$, while the second term equals to $0$.
Then we conclude the result.
\end{pf}

\par Moreover, we can determine the asymptotic behavior of the traces of $U_{\theta}$ on $\partial \Omega$
as follows.

\begin{lemma}
There holds the following
\begin{equation}
\sup_{\hat \lambda x\le y \le g(x)}|U_{\theta}(x-,y)-B_{\theta}(x-)|\to 0 \qquad  \text{as $x\to \infty$}
\end{equation}
for any $\hat \lambda \in(sup \lambda_1, \inf g')$.
\end{lemma}

\begin{pf}
Notice that
\begin{eqnarray*}
&&\sup_{\hat \lambda x\le y \le g(x)}|U_{\theta}(x-,y)-B_{\theta}(x-)|\\
&&\le \sup_{\hat \lambda x
\le y \le g(x)}|U_{\theta}(x-,y)-\tilde{U}_{\theta}(x-,y)|
+\sup_{\hat \lambda x\le y \le g(x)}|\tilde{U}_{\theta}(x-,y)-B_{\theta}(x-)|.
\end{eqnarray*}

By Lemma 7.3, the first term on the right-hand side tends to zero.
In the same way as in \cite{zhang2}, the second term also tends to zero.
The proof is completed.
\end{pf}

\par From Lemmas 7.3 and 7.6, it follows that

\begin{lemma}
Let
\begin{equation}
B_{\theta}(\infty)=\lim_{x\to \infty}B_{\theta}(x-)
\end{equation}
and let
\begin{equation}
g'(\infty)=\lim_{x\to \infty}g'_{+}(x).
\end{equation}
Then
\begin{equation}
\lim _{x\to \infty}\sup_{\hat \lambda x\le y \le g(x)}|\lambda_1(U_{\theta}(x-,y))-\lambda_1(B_{\theta}(x-))|= 0,
\end{equation}
and
\begin{equation*}
B_{\theta}(\infty)\cdot (-g'(\infty),1)=0.
\end{equation*}
\end{lemma}

\par Repeating the argument as in \cite{liu2} and by Lemmas 7.3 and 7.7, we can prove

\begin{lemma}
\par Let $U_{\infty}=\lim_{y\to -\infty}U_0(y)$ for the initial data $U_0(y)$ at $x=0$.
\begin{enumerate}
\item[\rm (i)] If $\lambda_1(B_{\theta}(\infty))>\lambda_1(U_{\infty})$, then
\begin{equation}
B_{\theta}(\infty)\in R_1^+(U_{\infty}).
\end{equation}

\item[\rm (ii)]
If $\lambda_1(B_{\theta}(\infty))\le \lambda_1(U_{\infty})$, then
\begin{equation}
B_{\theta}(\infty)\in S_1^-(U_{\infty}).
\end{equation}
\end{enumerate}
Therefore, the equation
\begin{equation}
\Phi(0,0,0,0,\alpha_{\infty}; U_{\infty})=B_{\theta}(\infty)
\end{equation}
has a unique solution $\alpha_{\infty}$.\label{ff}
\end{lemma}

\par Considering the geometry of the boundary and performing
the same way as in \cite{zhang2}, we can obtain
\begin{lemma}
Suppose that $|g'(\infty)|$ is small, then
\begin{enumerate}
\item[\rm (i)] If $g'(\infty)<0$, then $\lambda_1(B_{\theta}(\infty))>\lambda_1(U_{\infty})$;
\item[\rm (ii)] If $g'(\infty)=0$, then $\lambda_1(B_{\theta}(\infty))=\lambda_1(U_{\infty})$;
\item[\rm (iii)] If $g'(\infty)>0$, then $\lambda_1(B_{\theta}(\infty))<\lambda_1(U_{\infty})$.
\end{enumerate}
\end{lemma}

\par By carrying out the same arguments as in \cite{zhang2} and employing the above lemmas,
we finally have the asymptotic behavior of entropy solutions.

\begin{theorem}
Suppose that $\mathrm{TV}(U_0)+\mathrm{TV}(g')$ is sufficiently small.
\begin{enumerate}
\item[\rm (i)] If $g'(\infty)<0$, then there exists a $1$-shock which approaches to
the shock wave with $(\alpha_{\infty},0,0,0,0)$ both in strength and speed as $x\to \infty$;
moreover, the total variation of $U_{\theta}$ outside this shock approaches to zero as $x\to \infty$.
\item[\rm (ii)] If $g'(\infty)=0$, then $\sup_{y<g(x)}|U_{\theta}(x,y)-U_{\infty}| \to 0$ as $x\to \infty$.
\item[\rm (iii)] If $g'(\infty)>0$, then the amount of shocks approaches to zero as $x\to \infty$
and $U(x,y)$ approaches the rarefaction wave with $(\alpha_{\infty},0,0,0,0)$, where
$(\alpha_{\infty},0,0,0,0)$ is given in Lemma {\rm \ref{ff}}.
\end{enumerate}
\end{theorem}

\section{\textbf{Supersonic Reacting Euler Flow past Lipschitz Wedge with Large Angle}}

Now we consider the general case when the wedge angle is arbitrary large, but less than the sonic angle. We establish a theory of global existence and asymptotic behavior of entropy solutions for the initial-boundary
value problem \eqref{dd3}--\eqref{dd4} for system \eqref{d1}--\eqref{d2} for which $v_0(-\infty)$ is not zero in general.

\subsection{\textbf{Initial-boundary value problem involving a strong shock}}
For the wedge with large vertex angle, as in \cite{chen2}, we choose a suitable coordinate system (by rotation when it is necessary) such that the wedge has the lower boundary $\{y=g(x),x\ge 0\}$ with
\begin{equation}
g(0)=g'(0)=0, \qquad g\in C([0,\infty]), \qquad g'\in \mathrm{BV}.
\end{equation}

\begin{center}
\setlength{\unitlength}{1mm}
\begin{picture}(120,70)(-2,-10)
\linethickness{1pt}
\put(30,0){\vector(0,1){50}}
\put(0,30){\vector(1,0){120}}
\put(3,3){\vector(2,1){20}}
\put(3,7){\vector(2,1){20}}
\put(3,11){\vector(2,1){20}}
\put(3,-1){\vector(2,1){20}}
\qbezier(30,30)(40,31)(50,33)
\qbezier(50,33)(90,25)(110,32)
\qbezier(30,30)(50,15)(90,2)
\put(90,4){$Shock$}
\put(24,26){$O$}
\put(116,26){$x$}
\put(27,50){$y$}
\put(-6,20){$(u_0(y), v_0(y))$}
\put(34,30){\line(1,2){2}}
\put(38,31){\line(1,2){2}}
\put(42,31){\line(1,2){2}}
\put(46,32){\line(1,2){2}}
\put(50,32){\line(1,2){2}}
\put(54,32){\line(1,2){2}}
\put(58,31){\line(1,2){2}}
\put(62,31){\line(1,2){2}}
\put(66,31){\line(1,2){2}}
\put(70,30){\line(1,2){2}}
\put(74,30){\line(1,2){2}}
\put(78,29){\line(1,2){2}}
\put(82,29){\line(1,2){2}}
\put(86,29){\line(1,2){2}}
\put(90,29){\line(1,2){2}}
\put(94,29){\line(1,2){2}}
\put(98,30){\line(1,2){2}}
\put(100,20){\vector(-2,1){16}}
\put(100,18){$y=g(x)$}
\put(52,10){$\Omega$}
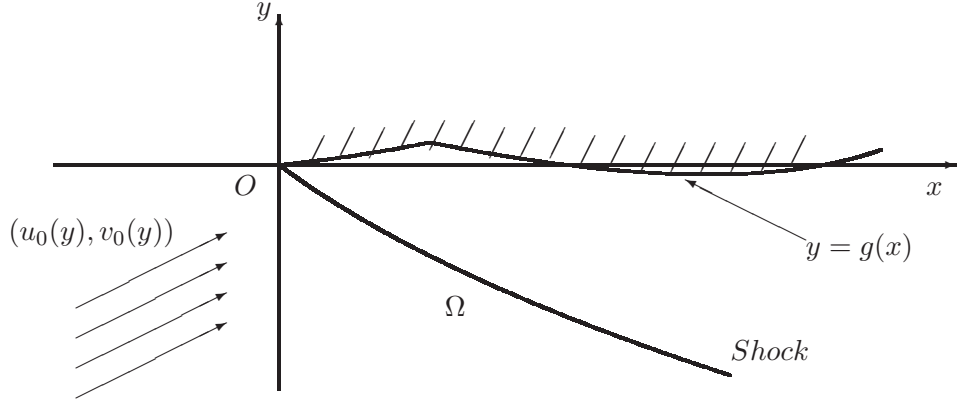
\put(18,-8){Fig. 11. Initial-boundary problem with large vertex angle}
\end{picture}
\end{center}

For the non-reaction problem with straight boundary $\{x\ge 0, y\equiv 0\}$ and uniform incoming flow $U_0(-\infty)$, if we assume that
\begin{equation}
0<\arctan\Big(\frac{v_0(-\infty)}{u_0(-\infty)}\Big)<\omega_{crit},
\end{equation}
then there exists a supersonic state $U_{+}=(u_{+},0,p_{+},\rho_{+},Z_{+})\in S_1(U_0(-\infty))$ with entropy condition $u_{+}<u_0(-\infty)$ such that
the corresponding non-reaction problem \eqref{5}--\eqref{6} has a shock solution with
a leading shock front issuing from the vertex (see Fig. 12).
\bigskip

\begin{center}
\setlength{\unitlength}{1mm}
\begin{picture}(90,60)(-2,-10)
\linethickness{1pt}
\put(30,0){\vector(0,1){50}}
\put(0,30){\vector(1,0){70}}
\put(3,3){\vector(1,1){20}}
\put(30,30){\line(1,-1){25}}
\put(50,20){\vector(1,0){20}}
\put(30,30){\line(1,2){2}}
\put(34,30){\line(1,2){2}}
\put(38,30){\line(1,2){2}}
\put(42,30){\line(1,2){2}}
\put(46,30){\line(1,2){2}}
\put(50,30){\line(1,2){2}}
\put(54,30){\line(1,2){2}}
\put(58,30){\line(1,2){2}}
\put(62,30){\line(1,2){2}}
\put(22,26){$O$}
\put(66,26){$x$}
\put(23,50){$y$}
\put(-3,20){$(u_0(-\infty),v_0(-\infty))$}
\put(54,14){$(u_{+},v_{+})$}
\put(52,2){Shock}
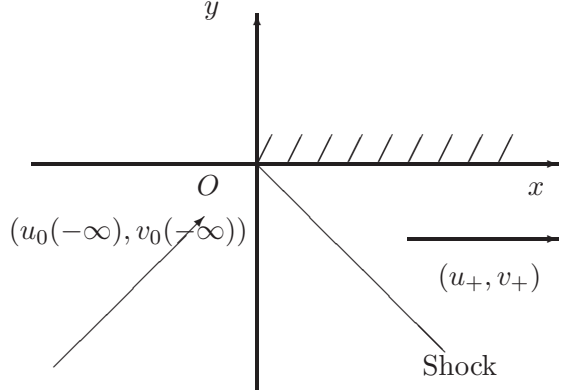
\put(-2,-8){Fig. 12. The background solution for the no-reaction problem}
\end{picture}
\end{center}
Moreover, there exist $r_1>0$ and $r_2>0$ such that, for any $U_1\in O_{r_2}(U_0(-\infty))$, the shock polar $S_1(U_1)\cap O_{r_1}(U_+)$ can be parameterized by the form
\begin{equation}
U=D(s,U_1) \qquad \text{with $U_+=D(s,U_{-\infty})$},
\end{equation}
where $s$ is the shock speed.

\subsection{\textbf{Riemann problem with a strong shock}}
To construct the approximate solutions, we need to solve the Riemann problem with a strong shock.

\begin{lemma}
Let $U_1\in O_{r_1}(U_0(-\infty))$ and $U_2\in O_{r_2}(U_+)$ with small positive constants $r_1>0$ and $r_2>0$. Then the Riemann problem
\begin{eqnarray}
\left\{
\begin{array}{ccc}
W(U)_{x}+H(U)_y=0,\\
 U|_{x=0}=
\begin{cases} U_1 &\quad y< y_0,\\U_2 &\quad y> y_0,
\end{cases}
\end{array} \right.
\label{rm2}
\end{eqnarray}
has a unique solution constituted by weak waves $\alpha_2,\alpha_3,\alpha_4,\alpha_5$, and a strong shock $s$, that is,
\begin{equation}
\Psi(\alpha_5,\alpha_4,\alpha_3,\alpha_2,0;D(s,U_1))=U_2.
\label{rm1}
\end{equation}
\end{lemma}

This lemma can be proved in the same way as in \cite{chen2} by solving (\ref{rm1}). Besides the Riemann problem for the interacting weak waves and the fractional steps in the previous sections, we also employ (\ref{rm2}) for dealing with the interaction between the weak waves and the strong wave. More precisely, we have the following lemma to include the strong shock.

\begin{center}
\setlength{\unitlength}{1mm}
\begin{picture}(50,90)(-2,-22)
\linethickness{1pt}
\multiput(0,-12)(0,3){26}{\line(0,2){2}}
\multiput(24,-12)(0,3){25}{\line(0,1){2}}
\multiput(48,-13)(0,3){25}{\line(0,2){2}}
\put(0,1){\line(1,1){20}}
\put(0,1){\line(3,-2){19}}
\put(0,1){\line(5,3){20}}
\put(0,1){\line(5,-2){24}}
\put(0,62){\line(5,-2){24}}
\put(0,62){\line(2,-3){16}}
\put(0,33){\line(5,-2){24}}
\put(24,52){\line(4,1){24}}
\put(24,23){\line(4,1){24}}
\put(24,-8){\line(4,1){24}}
\put(24,-8){\line(2,3){18}}
\put(24,-8){\line(1,1){18}}
\put(24,-8){\line(5,-1){18}}
\put(0,62){\line(1,-1){18}}
\put(9,54){$U_a$}
\put(38,33){$U_a$}
\put(2,24){$U_m$}
\put(2,-6){$U_b$}
\put(25,-13){$U_b$}
\put(-16,-20){Fig. 13. Interaction with the strong wave below}
\put(38,2){$\delta_{2(3,4)}$}
\put(36,18){$\delta_{5}$}
\put(44,-13){$s'$}
\put(12,35){$\beta_1$}
\put(14,50){$\beta_{2(3,4)}$}
\put(12,18){$\alpha_{5}$}
\put(12,6){$\alpha_{2(3,4)}$}
\put(15,-13){$s$}
\end{picture}
\end{center}

\begin{lemma}
Suppose that $U_b\in O_{r_1}(U_0(-\infty))$ and $U_a$, $U_m\in O_{r_2}(U_+)$ with
\begin{eqnarray}
\{U_m,U_a\}=(\beta_1,\beta_2,\beta_3,\beta_4,0),\\
\{U_b,U_m\}=(s,\alpha_2,\alpha_3,\alpha_4,\alpha_5),
\end{eqnarray}
and
\begin{equation}
\{U_b,U_a\}=(s',\delta_2,\delta_3,\delta_4,\delta_5).
\end{equation}
Then
\begin{eqnarray}
&& s'=s+K_{s_1}\beta_1+O(1)\Delta,\\
&& \delta_j=\alpha_j+\beta_j+K_{s_j}\beta_1+O(1)\Delta, \qquad \text{$j=2,3,4$},\\
&& \delta_5=\alpha_5+K_{s_5}\beta_1+O(1)\Delta,
\end{eqnarray}
with
\begin{equation}
|K_{s5}|<1,\qquad \text{$\sum_j|K_{sj}|\le M$ $\,\,\,$ for some $M>0$},
\end{equation}
and
\begin{equation}
\Delta=|\alpha_5|(|\beta_2|+|\beta_3|+|\beta_4|).
\end{equation}
\label{lem7}
\end{lemma}

\begin{lemma}
Suppose that
\begin{eqnarray}
\{U_b,U_m\}=(\alpha_1,\alpha_2,\alpha_3,\alpha_4,\alpha_5),\qquad
\{U_m,U_a\}=(s,\beta_2,\beta_3,\beta_4,\beta_5),
\end{eqnarray}
and
\begin{equation}
\{U_b,U_a\}=(s',\delta_2,\delta_3,\delta_4,\delta_5),
\end{equation}
with $U_b$, $U_m\in O_{r_2}(U_0(-\infty))$ and $U_a\in O_{r_1}(U_{+})$.
Then
\begin{eqnarray*}
s'=s+K_{s_1}\alpha_1+O(1)\sum_{j=1}^{5}|\alpha_j|,\qquad
\delta_j=\beta_j+O(1)\sum_{j=1}^5|\beta_j|.
\end{eqnarray*}
\label{lem8}
\end{lemma}

\begin{pf}
Actually, if we set $\alpha_j=0$ for all $j$, then $s'=s$ and $\delta_j=\beta_j$ for all $j$. Then the result follows.
\end{pf}
\bigskip

\begin{center}
\setlength{\unitlength}{1mm}
\begin{picture}(50,74)(-2,-4)
\linethickness{1pt}
\multiput(0,0)(0,3){24}{\line(0,2){2}}
\multiput(24,0)(0,3){24}{\line(0,2){2}}
\multiput(48,0)(0,3){24}{\line(0,2){2}}
\put(0,54){\line(3,2){20}}
\put(0,54){\line(3,-1){20}}
\put(0,54){\line(5,2){20}}
\put(0,18){\line(1,1){20}}
\put(0,18){\line(5,3){20}}
\put(0,18){\line(4,-1){20}}
\put(24,40){\line(5,-1){20}}
\put(24,40){\line(3,2){20}}
\put(24,40){\line(5,1){20}}
\put(2,70){$U_a$}
\put(2,30){$U_m$}
\put(2,5){$U_b$}
\put(36,15){$U_b$}
\put(36,60){$U_a$}
\put(12,66){$\beta_5$}
\put(12,56){$\beta_{2(3,4)}$}
\put(12,46){$s$}
\put(12,36){$\alpha_5$}
\put(12,23){$\alpha_{2(3,4)}$}
\put(12,11){$\alpha_1$}
\put(40,54){$\delta_5$}
\put(38,46){$\delta_{2(3,4)}$}
\put(40,38){$s'$}
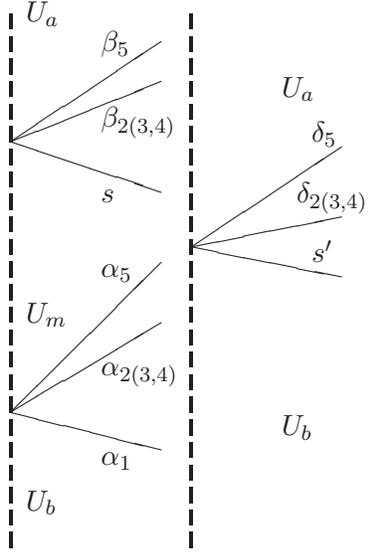
\put(-10,-6){Fig. 14. Interaction with the strong wave above}
\end{picture}
\end{center}

\medskip
\subsection{\textbf{Glimm-type functional involving the strong shock}}
We use the same grid points and mesh curves as in the previous sections. For the strip $\Omega_k$, we denote the strong shock in $\Omega_k$ by $s_k$. Without confusion, we also denote its speed by $s_k$ and its location by $y=\chi_k(x)$.
\par Let
\begin{equation}
\Omega_{k+}=\{\chi_k(x)<y\}\cap \Omega_k, \qquad \text{$\Omega_{k-}=\{\chi_k(x)>y\}\cap \Omega_k$}.
\end{equation}
For $J_k<J<J_{k+1}$, we denote $J_+=J\cap \Omega_{k+}$ and $J_-=J\cap \Omega_{k-}$.

\begin{definition}
\begin{eqnarray*}
&&L_j(J_{\pm})=\sum \{|\alpha|: \text{$\alpha$ is weak $j$-wave crossing $J_{\pm}$}\},\\[1.5mm]
&&Q(J_{\pm})=\sum \{|\alpha||\beta|: \text{$\alpha, \beta$ are weak waves, approaching and crossing $J_{\pm}$}\},\\[1.5mm]
&&L(J_+)=K_0^*L_0(J)+L_1(J_+)+K_2^*L_2(J_+)+K_3^*L_3(J_+)+K_4^*L_4(J_+)+K_5^*L_5(J_+),\\[1.5mm]
&&L(J_-)=L_1(J_-)+K_2^{**}L_2(J_-)+K_3^{**}L_3(J_-)+K_4^{**}L_4(J_-)+K_5^{**}L_5(J_-),\\[1.5mm]
&&F(J)=L(J_+)+KL(J_-)+K'Q(J_+)+KK^{''}Q(J_-),\\[1.5mm]
&&F_s(J)=|s_J-s_*|+C_*F(J),
\end{eqnarray*}
where $K, K', K^{''}, K_j^*, K_j^{**}$, and $C_*$ are all positive constants with
$$
K_0^*>|K_{b0}|, \qquad |K_{b5}|<K_5^*<\frac{1}{|K_{s_5}|}.
$$
\end{definition}

\begin{proposition}
Let $J_k<I<J<\tilde{J}_{k+1}$ such that $J$ is an immediate successor of $I$. Suppose that
\begin{eqnarray*}
&&\big{|}s_I-s_*\big{|}<\varepsilon,\\[1.5mm]
&&\big{|}U_{h,\theta}|_{I_+}-U_+\big{|}<\varepsilon_1,\\[1.5mm]
&&\big{|}U_{h,\theta}|_{I_-}-U_0(-\infty)\big{|}<\varepsilon_2
\end{eqnarray*}
for some $\varepsilon, \varepsilon_1$, and $\varepsilon_2>0$.
Then there exist positive constants $K, K',K^{''}, K_j^*, K_j^{**}, C_*$,
and $\tilde{\varepsilon}$, which are independent of $I, J$, and $k$,
such that, if $F_s(I)<\tilde{\varepsilon}$, then
\begin{equation*}
F_s(J)<F_s(I).
\end{equation*}
Furthermore, we have
\begin{eqnarray*}
&&\big{|}s_J-s_*\big{|}<\varepsilon,\\[1.5mm]
&&\big{|}U_{h,\theta}|_{J_+}-U_+\big{|}<\varepsilon_1,\\[1.5mm]
&&\big{|}U_{h,\theta}|_{J_-}-U_0(-\infty)\big{|}<\varepsilon_2.
\end{eqnarray*}
\end{proposition}

\begin{pf}
We consider only the case near the strong $1-shock$, since the other cases can be treated in the same way as in the previous sections.

\par Let $\Lambda$ be the diamond domain between the mesh curves $I$ and $J$.

\medskip
\par {\textbf{Case 1}}:
By Lemma \ref{lem7}, we have
\begin{eqnarray*}
&& L_1(J_+)-L_1(I_+)=-|\beta_1|,\\[1.5mm]
&& L_j(J_+)-L_j(I_+)\le |K_{s_j}||\beta_1|+O(1)\Delta, \qquad \text{$j=2,3,4$},\\[1.5mm]
&& L_5(J_+)-L_5(I_+)\le |K_{s_5}||\beta_1|+O(1)\Delta, \\[1.5mm]
&& L(J_-)-L(I_-)=0,\\[1.5mm]
&& Q(J_-)-Q(I_-)=0.
\end{eqnarray*}

\begin{center}
\setlength{\unitlength}{1mm}
\begin{picture}(50,92)(-2,-22)
\linethickness{1pt}
\multiput(0,-12)(0,3){26}{\line(0,2){2}}
\multiput(24,-12)(0,3){25}{\line(0,1){2}}
\multiput(48,-13)(0,3){25}{\line(0,2){2}}
\put(0,1){\line(1,1){20}}
\put(0,1){\line(3,-2){19}}
\put(0,1){\line(5,3){20}}
\put(0,1){\line(5,-2){24}}
\put(0,62){\line(5,-2){24}}
\put(0,62){\line(2,-3){16}}
\put(0,33){\line(5,-2){24}}
\put(24,52){\line(4,1){24}}
\put(24,23){\line(4,1){24}}
\put(24,-8){\line(4,1){24}}
\put(24,-8){\line(2,3){18}}
\put(24,-8){\line(1,1){18}}
\put(24,-8){\line(5,-1){18}}
\put(0,62){\line(1,-1){18}}
\put(5,-20){Fig. 15. Case 1}
\put(38,2){$\delta_{2(3,4)}$}
\put(36,18){$\delta_{5}$}
\put(38,-15){$s_{k+1}$}
\put(12,35){$\beta_1$}
\put(14,50){$\beta_{2(3,4)}$}
\put(12,18){$\alpha_{5}$}
\put(12,6){$\alpha_{2(3,4)}$}
\put(15,-13){$s_k$}
\end{picture}
\end{center}
Then we conclude that
\begin{equation*}
L(J_+)-L(I_+)\le (-1+\sum_{j=2}^5K_j^*|K_{s_j}|)|\beta_1|+O(1)\Delta,
\end{equation*}
and
\begin{equation*}
Q(J_+)-Q_(I_+)\le O(1)\Delta+O(1)|\beta_1|L_1(J_+).
\end{equation*}
Moreover,
\begin{equation*}
|s_{k+1}-s_*|\le |s_k-s_*|+|K_{s_1}||\beta_1|+O(1)\Delta.
\end{equation*}
Combining this with the above estimates, and choosing suitable constants $K_j^*$ and large constants $K$, $K^{\prime}$, and $K^{\prime\prime}$,
we conclude
\begin{equation*}
F_s(J)\le F_s(I), \quad \text{for $F_s(I)\le \tilde{\varepsilon}$}.
\end{equation*}

\par {\textbf{Case 2}}: By Lemma \ref{lem8}, we have
\begin{eqnarray*}
&& s_{k+1}=s_k+O(1)|\boldsymbol{\beta}|, \\[1.5mm]
&& \delta_j=\alpha_j+O(1)|\boldsymbol{\beta}|,  \qquad j=1, \cdots, 5,
\end{eqnarray*}
where $|\boldsymbol{\beta}|=\sum_{j=1}^5|\beta_j|$.
Then
\begin{equation*}
L(J_-)-L(I_-)\le -|\boldsymbol{\beta}|
\end{equation*}
for suitable choice of constants $K_j^{**}$. By choosing sufficiently large $K$, we finally have the desired result.

The proof is complete.
\bigskip

\begin{center}
\setlength{\unitlength}{1mm}
\begin{picture}(50,74)(-2,-4)
\linethickness{1pt}
\multiput(0,0)(0,3){24}{\line(0,2){2}}
\multiput(24,0)(0,3){24}{\line(0,2){2}}
\multiput(48,0)(0,3){24}{\line(0,2){2}}
\put(0,54){\line(3,2){20}}
\put(0,54){\line(3,-1){20}}
\put(0,54){\line(5,2){20}}
\put(0,18){\line(1,1){20}}
\put(0,18){\line(5,3){20}}
\put(0,18){\line(4,-1){20}}
\put(24,40){\line(5,-1){20}}
\put(24,40){\line(3,2){20}}
\put(24,40){\line(5,1){20}}
\put(12,66){$\beta_5$}
\put(12,56){$\beta_{2(3,4)}$}
\put(12,46){$s_k$}
\put(12,36){$\alpha_5$}
\put(12,23){$\alpha_{2(3,4)}$}
\put(12,11){$\alpha_1$}
\put(40,54){$\delta_5$}
\put(38,46){$\delta_{2(3,4)}$}
\put(40,38){$s_{k+1}$}
\put(10,-6){Fig. 16. Case 2}
\end{picture}
\end{center}

\end{pf}

\medskip
\subsection{\textbf{Estimates of reaction steps for the strong shock}}
By Lemma \ref{q1}, we have
\begin{equation*}
U_b-\tilde{U}_b=\|Z_0\|_{\infty}e^{-\Phi_1 kh}O(h),
\end{equation*}
\begin{equation*}
U_a-\tilde{U}_a=\|Z_0\|_{\infty}e^{-\Phi_1 kh}O(h).
\end{equation*}
Then
\begin{equation}
\tilde{s}_k-s_k=\|Z_0\|_{\infty}e^{-\Phi_1 kh}O(h).
\end{equation}

\medskip
\begin{center}
\setlength{\unitlength}{1mm}
\begin{picture}(100,52)(-2,-10)
\linethickness{1pt}
\multiput(0,0)(0,3){13}{\line(0,2){2}}
\multiput(25,0)(0,3){13}{\line(0,2){2}}
\multiput(73,-1)(0,3){13}{\line(0,2){2}}
\multiput(98,-1)(0,3){13}{\line(0,2){2}}
\put(73,18){\line(4,-1){20}}
\put(0,18){\line(4,-1){20}}
\put(19,9){$\tilde{s}_k$}
\put(14,20){$\tilde{U}_a$}
\put(92,10){$s_k$}
\put(87,20){$U_a$}
\put(10,6){$\tilde{U_b}$}
\put(83,6){$U_b$}
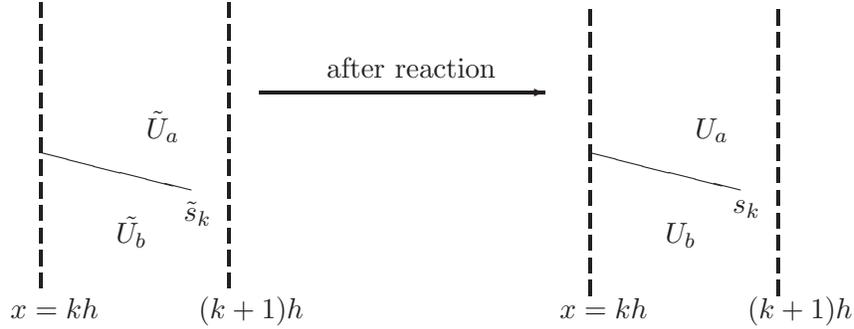
\put(8,-11){Fig. 17.  Change of the strength of the strong shock}
\put(38,28){after reaction}
\put(29,26){\vector(1,0){38}}
\put(-4,-4){$x=kh$}
\put(69,-4){$x=kh$}
\put(21,-4){$(k+1)h$}
\put(94,-4){$(k+1)h$}
\end{picture}
\end{center}

\bigskip
\medskip
As in the previous sections, we still have
\begin{equation}
 F_s(J_k)-F_s(\tilde{J_k})\le Ch\|w_{5,0}\|_{\infty}e^{-\Phi_1 kh}(F_s(\tilde{J_k})+1)^2.
 \end{equation}
This gives the uniform bounds on $F_s(J_k)$.

\medskip
\subsection{\textbf{Global existence and asymptotic behavior of entropy solutions for
the Lipschitz wedge with large angle}}

\medskip
We finally have the following theorem.

\begin{theorem}
Suppose that $0<\arctan\big(\frac{v_0(-\infty)}{u_0(-\infty)}\big)<\omega_{crit}$. If $\mathrm{TV}(W(U_0))+\mathrm{TV}(g')$ is sufficiently small,
then the fractional-step Glimm scheme can generate
a family of approximate solutions $U_{h,\theta}(x,y)$ that have uniformly bounded variation in the $y$--direction. Moreover, there exists a null set $N\subset \Pi_{k=0}^{\infty}(-1,1)$ such that, for every $\theta \in \Pi_{k=0}^{\infty}(-1,1)\setminus N$, there exist a sequence $h_j\to 0$ such that
\begin{equation}
U_{\theta}\stackrel{L_{loc}^1}{=}\lim _{h_i\to 0}U_{h_i,\theta}
\end{equation}
is a weak solution to problem \eqref{dd3}--\eqref{dd4} for system \eqref{d1}--\eqref{d2}.
Moreover,  $U_{\theta}$ has uniformly  bounded variation in the $y$--direction.
\end{theorem}

\par In the same way as in \cite{chen2}, we have
\begin{theorem}[Asymptotic behavior]
Let $\omega_{\infty}=\lim_{x\to \infty} \arctan(g'(x+))$. Then
\begin{equation}
\lim_{x\to \infty}\sup_{\chi_{\theta}(x)<y<g(x)}\big|
\arctan\big(\frac{v_{\theta}(x,y)}{u_{\theta}(x,y)}-\omega_{\infty}\big)\big|=0,
\end{equation}
and
\begin{equation}
\lim_{x\to \infty}\sup_{y<\chi_{\theta}(x)}
\big|\arctan\big(\frac{v_{\theta}(x,y)}{u_{\theta}(x,y)}\big)\big|=0.
\end{equation}
\end{theorem}

\bigskip
\bigskip
\medskip
\noindent
\textbf{Acknowledgements:}
The research of
Gui-Qiang Chen was supported in part by the UK EPSRC Science and Innovation
Award to the Oxford Centre for Nonlinear PDE (EP/E035027/1),
the NSFC under a joint project Grant 10728101, and
the Royal Society--Wolfson Research Merit Award (UK).
Changguo Xiao was supported in part by  the NSFC under a joint project Grant 10728101.
Yongqian Zhang was supported in part by NSFC Project 11031001, NSFC Project 11121101,
and the 111 Project B08018 (China).

\bigskip
\bigskip
\bigskip
\bigskip
\newpage

\centerline{\bf References}

\bibliographystyle{plain}

\begin{thebibliography}{99}

\bibitem{chen1} G.-Q. \textsc{Chen} \textsc{and} D. \textsc{Wagner},
  \emph{Global entropy solutions to exothermically reacting, compressible Euler equations},
  J. Differential Equations, {\bf 191} (2003), 277--322.

\bibitem{chen2} G.-Q. \textsc{Chen}, Y. Q. \textsc{Zhang}, \textsc{and} D. W. \textsc{Zhu},
 \emph{Existence and stability of supersonic Euler flows past Lipschitz wedges},
 Arch. Rational Mech. Anal. {\bf 181} (2006), 261--310.

\bibitem{chen3} G.-Q. \textsc{Chen}, Y. Q. \textsc{Zhang}, \textsc{and} D. W. \textsc{Zhu},
\emph{Stabilily of compressible vortex sheets in steady supersonic Euler flows over Lipschitz walls},
SIAM J. Math. Anal. {\bf 38} (2007), 1660--1693.

\bibitem{sxchen} S.-X. \textsc{Chen}, \emph{Asymptotic behavior of supersonic flow
past a convex combined wedge}, Chin. Ann. Math. {\bf 19B:3} (1998), 255--264.

\bibitem{courant} R. \textsc{Courant} \textsc{and} K. O. \textsc{Friedrichs},
\emph{Supersonic Flow and Shock Waves}, Wiley-Interscience, New York, 1948.

\bibitem{dafermos1} C. \textsc{Dafermos}, \emph {Hyperbolic Conservation Laws
in Continuum Physics}, Springer-Verlag: Berlin, 2005.

\bibitem{dafermos2} C. \textsc{Dafermos} \textsc{and} \textsc{L.}  \textsc{Hsiao},
\emph {Hyperbolic systems of balance laws with inhomegeneity and dissipation},
Indiana Univ. Math. J. {\bf 31} (1982), 471--491.

\bibitem{glimm} J. \textsc{Glimm},
\emph{Solutions in the large for nonlinear hyperbolic systems of equations},
Comm. Pure Appl. Math. {\bf 18} (1965), 697--715.

\bibitem{lax} P. D.  \textsc{Lax},
\emph{Hyperbolic systems of conservation laws \uppercase\expandafter{\romannumeral 2}},
Comm. Pure Appl. Math. {\bf 10} (1957), 537--566.

\bibitem{liu1} T.-P.  \textsc{Liu},
\emph{Solutions in the large for the equations of nonisentropic gas dynamics},
Indiana Univ. Math. J. {\bf 26} (1977), 147--177.

\bibitem{liu2} T.-P.  \textsc{Liu},
\emph{Large-time behaviour of initial and initial-boundary value problems of
a general systems of hyperbolic conservation laws},
Comm. Math. Phys. {\bf 55} (1977), 163--177.

\bibitem{luskin} M. \textsc{Luskin} \textsc{and} J. B. \textsc{Temple},
\emph{The existence of a global weak solution to the nonlinear waterhammer problem},
Comm. Pure Appl. Math. {\bf 34} (1982), 697--735.

\bibitem{smoller1} J. \textsc{Smoller}, \emph{Shock Waves and Reaction-Diffusion Equations},
Springer-Verlag: New York, 1983.

\bibitem{temple1} J. B. \textsc{Temple},
\emph{Solutions in the large for the nonlinear hyperbolic conservation laws of gas dynamics},
J. Differential Equations, {\bf 41} (1981), 96--161.

\bibitem{volpert} A. I. \textsc{Volpert},
\emph{The space BV and quasilinear equations},
 Mat. Sb. (N.S), {\bf 73} (1967), 255--302 (in Russian);
 Math. USSR Sb. {\bf 2} (1967), 225--267 (in English).

\bibitem{ying1} L.-A.  \textsc{Ying} \textsc{and} C.-H. \textsc{Wang},
\emph{Global solutions of the Cauchy problem for a nonhomogeneous quasilinear hyperbolic system},
Comm. Pure Appl. Math. {\bf 33} (1980), 579--597.

\bibitem{ying2} L.-A.  \textsc{Ying} \textsc{and} C.-H. \textsc{Wang},
\emph{Solutions in the large for nonhomogeneous quasilinear hyperbolic systems of equations},
J. Math. Anal. Appl. {\bf 78} (1980), 440--454.

\bibitem{zhang1} Y. Q. \textsc{Zhang},
\emph{Global existence of steady supersonic potential flow past a curved wedge with piecewise smooth boundary},
SIAM J. Math. Anal. {\bf 31} (1999), 166--183.

\bibitem{zhang2} Y. Q. \textsc{Zhang},
\emph{Steady supersonic flow past an almost straight wedge with large vertex angle},
J. Differential Equations, {\bf  192} (2003), 1--46.
\end{thebibliography}

\end{document}